\documentclass[10pt, journal]{IEEEtran}

\usepackage{cmap}

\usepackage{graphicx}      

\makeatletter
\let\NAT@parse\undefined
\makeatother

\usepackage[numbers]{natbib}

\usepackage[T1]{fontenc}
\usepackage[utf8]{inputenc}

\usepackage{caption}

\usepackage{subcaption}

\usepackage{float}

\usepackage{mathtools} 
\mathtoolsset{showonlyrefs}

\usepackage{breqn}

\usepackage{afterpage}

\usepackage{amsmath}%

\usepackage{enumitem}

\usepackage{multirow}
\usepackage{epstopdf}
\usepackage{url}

\usepackage{xcolor}

\usepackage{tikz}
\usetikzlibrary{arrows.meta,shapes.arrows}
\usetikzlibrary{chains,shapes.multipart}
\usepackage{tkz-graph}
\usetikzlibrary{arrows, shapes}
\usetikzlibrary{external,calc,patterns,decorations.pathmorphing,decorations.markings}
\usetikzlibrary{external}

\usepackage{tabulary}
\usepackage{booktabs}

\usepackage{comment}

\usepackage{amsfonts}
\DeclareMathSymbol{\shortminus}{\mathbin}{AMSa}{"39}

\usepackage{amsthm}

\theoremstyle{plain}
\newtheorem{prop}{Proposition}

\definecolor{vehblue}{RGB}{25,105,188}
\definecolor{vehbrown}{RGB}{127,51,0}
\definecolor{vehpink}{RGB}{178,0,255}

\graphicspath{ {C:/Images-KTH/} {figures/} } 
\begin{document}



\title{\LARGE \bf
Coordinating Vehicle Platoons for \\ Highway Bottleneck Decongestion and Throughput Improvement
}

\author{Mladen \v{C}i\v{c}i\'{c}$^1$, \and Li Jin$^2$, \and Karl Henrik Johansson$^1$
\thanks{
$\!\!\!\!{}^1$ Division of Decision and Control Systems, 
School of Electrical Engineering and Computer Science, 
KTH Royal Institute of Technology, 
Malvinas v{\"a}g 10, 10044 Stockholm, Sweden, \newline
Mladen \v{C}i\v{c}i\'{c}: {\tt\small cicic@kth.se}, Karl Henrik Johansson: {\tt\small kallej@kth.se}
\newline
${}^2$ Department of Civil and Urban Engineering,
C2SMART University Transportation Center,
New York University Tandon School of Engineering,
15 MetroTech Center, Brooklyn, NY 11201, USA, \newline
Li Jin: {\tt\small lijin@nyu.edu}\newline
The research leading to these results has received funding from the European Union's Horizon 2020 research and innovation programme under the Marie Sk{\l}odowska-Curie grant agreement No 674875, NYU Tandon School of Engineering, the C2SMART University Transportation Center, VINNOVA within the FFI program under contract 2014-06200, the Swedish Research Council, the Swedish Foundation for Strategic Research, and Knut and Alice Wallenberg Foundation. The first and third authors are affiliated with the Wallenberg AI, Autonomous Systems and Software Program (WASP).}

}%
\renewcommand\footnotemark{}
\renewcommand\footnoterule{}

\maketitle

\newcommand{\shat}{\hat{\sigma}}
\newcommand{\bro}{\bar{\rho}}
\newcommand{\brr}{\bar{r}}
\newcommand{\barf}{\bar{\varphi}}
\newcommand{\K}{\mathcal{K}}

\begin{abstract}

%
%
%

Truck platooning is a technology that is expected to become widespread in the coming years.
Apart from the numerous benefits that it brings, its potential effects on the overall traffic situation need to be studied further, especially at bottlenecks and ramps.
Assuming we can control the platoons from the infrastructure, they can be used as controlled moving bottlenecks, actuating control actions on the rest of the traffic, and potentially improving the throughput of the whole system.
In this paper, we use a multi-class cell transmission model to capture the interaction between truck platoons and background traffic, and propose a corresponding queuing model, which we use for control design.
We use platoon speeds, and the number of lanes platoons occupy as control inputs, and design a control strategy for throughput improvement of a highway section with a bottleneck.
By postponing and shaping the inflow to the bottleneck, we are able to avoid traffic breakdown and capacity drop, which significantly reduces the total time spent of all vehicles.
We derived the estimated improvement in throughput that is achieved by applying the proposed control law, and then tested it in a simulation study and found that the median delay of all vehicles by $75.6\%$ compared to the uncontrolled case.
Notably, although they are slowed down while actuating control actions, platooned vehicles experience less delay compared to the case without control, since they avoid going through congestion at the bottleneck.
\end{abstract}


\section{Introduction}


With truck platooning progressing persistently towards becoming a commonplace technology \citep{bergenhem2012overview}, studying and understanding the impact it will have on the overall traffic is becoming increasingly important.
Apart from providing potential fuel savings through air drag reduction \citep{bonnet2000fuel}, which was traditionally seen as its primary purpose, as well as having the potential to greatly reduce the work load on drivers \citep{mckinsey2018platooning}, truck platooning is also expected have a positive impact on traffic efficiency through reducing the headways between vehicles \citep{shladover2012impacts, lioris2017platoons}, alleviating the adverse effect trucks have on the traffic \citep{moridpour2015impact}.
Although there have been numerous field tests of truck platooning in real traffic \citep{alam2010experimental, aarts2016european}, insufficient emphasis has been put on understanding how these platoons affect the behavior of other vehicles on the road; thus the possible drawbacks of this technology are not yet fully understood \citep{bhoopalam2018planning}.



One identified problem pertains to the interaction between truck platoons and passenger cars close to on- and off-ramps, and bottlenecks in general \citep{rijkswaterstaat2016platooning}.
There is concern that long platoons might block access to an off-ramp, or from an on-ramp, forcing drivers to slow down excessively or cut into a platoon, resulting in significant disturbances for both the platoon and the rest of the traffic.
Furthermore, the arrival of platoons can cause traffic breakdown at a bottleneck, causing reduction of throughput due to the capacity drop phenomenon.
Recently, there have been efforts to address this problem in microscopic \citep{duret2019hierarchical} and macroscopic \citep{jin2018modeling} frameworks.
In this paper, we are focusing on applying a new type of macroscopic control, using the truck platoons as actuators.


Bottleneck decongestion	has long been tackled by classical traffic control measures, such as ramp metering \citep{wang2014local} and variable speed limits \citep{hadiuzzaman2012variable}.
However, both of these control methods require additional fixed equipment to be installed upstream of the bottleneck, which limits their flexibility, especially when it comes to handling temporary bottlenecks, such as work zones, incidents etc.
With the introduction of connected autonomous vehicles to the highways, new opportunities for sensing \citep{herrera2010evaluation} and actuation \citep{cui2017stabilizing, wu2018stabilizing} of the traffic are becoming available.
Lagrangian actuation, where we use a subset of vehicles that can be controlled directly from the infrastructure to restrict the traffic flow, is lately garnering some attention \citep{piacentini2018traffic, vinitsky2018lagrangian, cicic2018traffic}.
This approach effectively emulates ramp metering and variable speed control, achieving a similar type of regulation without the need for additional fixed equipment, allowing us to also control areas away from known permanent bottlenecks.
Moving bottleneck control is one such method, where we use slower moving vehicles to restrict the mainstream traffic flow at some points, delaying the arrival of some vehicles in a way similar to ramp metering and variable speed limits.

Due to their large size and the existence of fleet management infrastructure, truck platoons are an ideal candidate for moving bottleneck control.
Since they consist of heavy, slow-moving vehicles, truck platoons will act as moving bottlenecks with or without external control, and we may use the communication channels already in place to send centrally computed reference speeds and other control actions \citep{alam2015heavy}.
This way, we are able to mitigate the negative effects trucks have on the traffic, and even improve the overall traffic situation.
Apart from these positive effects on the traffic, truck platoons may improve the situation for themselves as well, leading to potentially less delay, smoother speed profiles, as well as increased predictability.	

%
%
%
%

%
%
%
%
%
%
%







To this end, we need an appropriate model of the mutual influence that truck platoons and the rest of the traffic have on each other, that is both tractable and sufficiently rich.
Microscopic traffic models allow for a fairly straightforward representation of trucks and platoons \citep{liu2016modeling}, and PDE models offer a consistent way of introducing moving bottlenecks \citep{delle2014scalar}, but both are overly complex and detailed for link-level control synthesis.
The multi-class cell transmission model (CTM) \citep{cicic2019multiclass} presents a good balance of complexity and tractability, and will therefore be used as a simulation model in this work.
We further simplify this model using a queuing representation, and use this newly derived model for control design.

%
%
%
%


The problem that we are addressing in this paper is bottleneck decongestion using randomly arriving platoons as actuators, with their speed and the number of lanes they occupy as a control input.
The main contributions of this work are the queuing-based model for predicting the evolution of the traffic, and the control law that uses this prediction to improve the throughput.
The designed control law is tested on a road segment upstream of a lane drop bottleneck that has one on-ramp and one off-ramp, and shown to achieve a significant reduction in total time spent.
We conduct basic stability analysis of the controlled system, and derive estimates for the improved throughput.
The median delay of all vehicles is reduced by $75.6\%$ in case the proposed control is applied, compared to the case with no control.
Even though the platoons are slowed down while actuating control actions, they experience overall less delay compared to the case without control, since they avoid going through congestion at the bottleneck.








The paper is structured as follows.
In Section~\ref{sec:model}, we present the multi-class CTM and introduce its simplified queuing representation.
Then, in Section~\ref{sec:control}, we use the said simplified model to design control laws for improving the throughput of the road, and in Section~\ref{sec:analysis} give a stability analysis of the proposed control law, as well as estimates on achieved throughput.
Section~\ref{sec:sim} describes the simulation setup and results, and finally, in Section~\ref{sec:conclusion} we conclude and discuss the results.

%
%
%
%
%

\section{Model}
\label{sec:model}

In this section, we present the traffic models that will be used for analysis, simulation and control design.
The base model, multi-class CTM, is augmented to properly represent the behavior of platoons moving slower than the rest of traffic.
Since this model still has a high number of states and control inputs, we propose a simplified queuing model, that is consistent with the multi-class CTM, and use it for control design.
Control actions will be calculated using predictions based on this simplified model, and then applied to the more complex simulation model for evaluations.

\subsection{The multi-class CTM}

The simulation model that is used in this work is a multi-class extension of the well-known CTM \citep{daganzo1994cell}, and it is a variant of the model used in \citep{piacentini2019mcctm} and \citep{cicic2019multiclass}.
Let $\mathcal{K}$ be the set of vehicle classes.
The traffic density of vehicles of class ${\kappa \in \mathcal{K}}$ in cell $i$ at time $t$ will be expressed in terms of passenger car equivalents, and denoted $\rho_i^\kappa(t)$.
We allow each of the classes to have a distinct free flow speed $U_i^\kappa(t)$ in every cell, varying in time.
These free flow speeds must not be higher than the overall maximum free flow speed for the cell, $U_i^{\kappa}(t) \le V_i$.
In practice, we use $U_i^{\kappa}(t)$ to capture some richer behaviour not covered by the base model, like platoons stop-and-go waves, as well as to apply the control action to the classes of vehicles we have control over.
We will use the platoon model given in \citep{cicic2019multiclass} for simulation and control design.	



\begin{figure}[b]
    \centering
    
        \begin{tikzpicture}
    \node[anchor=south west,inner sep=0] (image) at (0,0) {\includegraphics[width=\columnwidth]{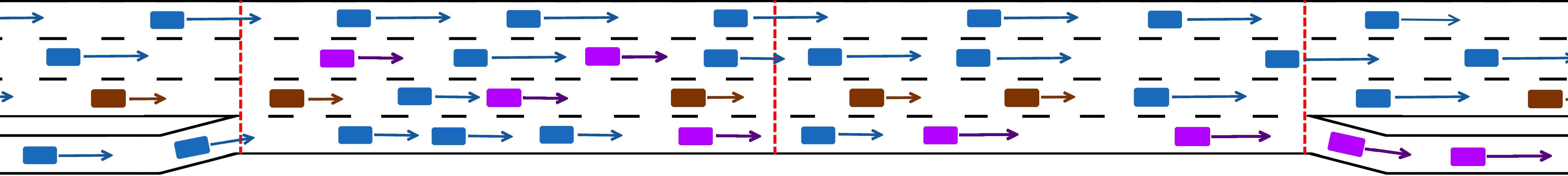}};
    \begin{scope}[x={(image.south east)},y={(image.north west)}]
        \draw[->] (0.11,1.1) -- +(0.1,0) node[above] {\scriptsize \hspace*{-0.7cm} \textcolor{vehbrown}{$\!\!q_{i\!-\!1}^a\!(t)$}, \textcolor{vehblue}{$\!\!q_{i\!-\!1}^b\!(t)$}, \textcolor{vehpink}{$\!\!q_{i\!-\!1}^c\!(t)$}};
        \draw[->] (0.44,1.1) -- +(0.1,0) node[above] {\scriptsize \hspace*{-0.7cm} \textcolor{vehbrown}{$\!\!q_{i}^a\!(t)$}, \textcolor{vehblue}{$\!\!q_{i}^b\!(t)$}, \textcolor{vehpink}{$\!\!q_{i}^c\!(t)$}};
        \draw[->] (0.78,1.1) -- +(0.1,0) node[above] {\scriptsize \hspace*{-1cm} \textcolor{vehbrown}{$\!\!q_{i\!+\!1}^a\!(t)$}, \textcolor{vehblue}{$\!\!q_{i\!+\!1}^b\!(t)$}};
        \node at (0.35,-0.05) {\scriptsize \textcolor{vehbrown}{$\!\!\rho_{i}^a\!(t)$}, \textcolor{vehblue}{$\!\!\rho_{i}^b\!(t)$}, \textcolor{vehpink}{$\!\!\rho_{i}^c\!(t)$}};
        \node at (0.65,-0.05) {\scriptsize \textcolor{vehbrown}{$\!\!\rho_{i\!+\!1}^a\!(t)$}, \textcolor{vehblue}{$\!\!\rho_{i\!+\!1}^b\!(t)$}, \textcolor{vehpink}{$\!\!\rho_{i\!+\!1}^c\!(t)$}};
        \draw[->] (0.1,-0.1) -- +(0.1,0.1) node[below] {\scriptsize \hspace*{-0.7cm} \textcolor{vehblue}{$\!\!r_{i}^b\!(t)$}};
        \draw[->] (0.79,0) -- +(0.1,-0.1) node[below = -0.1cm] { \scriptsize \hspace*{-0.9cm}  \textcolor{vehpink}{$\!\!s_{i\!+\!1}^c\!(t)$}};
    \end{scope}
\end{tikzpicture}    
    \caption{An example of three-class traffic flows in two cells. Vehicle classes $a$, $b$, $c$ are colour-coded. Cell $i$ receives traffic of all three classes from cell $i-1$ and class $b$ vehicles from an on-ramp. Class $a$ and $b$ vehicles are downstream-bound, and will leave cell $i+1$ and enter cell $i+2$, whereas class $c$ vehicles are off-ramp-bound and will leave cell $i+1$ via the off-ramp.}
    \label{fig:twocell}
\end{figure}
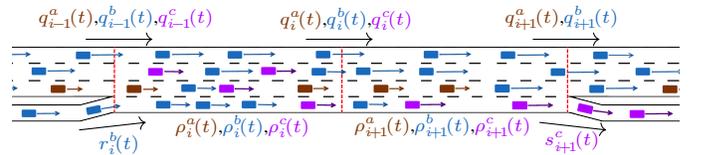

Consider a highway stretch consisting of $N$ cells.
The evolution of cell traffic densities for each class is given by
\begingroup
\medmuskip=0.3\medmuskip
\thickmuskip=0.3\thickmuskip
\begin{equation}
\label{eq:MCTM}
\rho_i^\kappa (t+1) = \rho_i^\kappa (t) + \frac{T}{L_i}\left(q_{i-1}^\kappa(t) - q_{i}^\kappa(t) + r_i^{\kappa}(t) - s_i^{\kappa}(t)\right),
\end{equation}
\endgroup
where $r_i^{\kappa}(t)$ is the inflow and $s_i^{\kappa}(t)$ the outflow of each vehicle class from a potential on-ramp and to a potential off-ramp, respectively.
An example of such traffic flows is given in Figure~\ref{fig:twocell}.
The traffic flow of each class from cell $i$ to cell $i+1$ is given by
\begin{equation}
q_i^\kappa(t) = \min\left\{D_i^\kappa(t), S_{i+1}^\kappa(t)\right\}.
\end{equation}

The demand and supply functions of each class $D_i^\kappa(t)$ and $S_i^\kappa(t)$ also depend on vehicles of other classes.
Denoting
\begin{align}
d_i^{\kappa}(t) &= U_i^{\kappa}(t) \rho_i^{\kappa}(t),\\
d_i(t) &= \sum\limits_{\kappa \in \mathcal{K}} d_i^{\kappa}(t), 
\end{align}
we write the demand and supply functions
\begin{align}
\label{eq:DMC}
D_i^\kappa(t) &= d_i^{\kappa}(t)\min\left\{1,\frac{Q_i(t)}{d_i(t)}\right\},\\
\label{eq:SMC}
S_{i}^\kappa(t) &= \frac{\rho_{i-1}^\kappa(t)}{\rho_{i-1}(t)} \min\left\{ W_i(P_{i} - \rho_{i}(t)), Q_i(t), F_{i-1}(t) \right\}.
\end{align}
Here, cell parameters $L_i$, $V_i$, $W_i$, $\sigma_i$ and $P_i$ are the length, free flow speed, congestion wave speed, critical density and jam density of cell $i$, respectively, and the cell capacity is given by
\begin{align}
Q_i(t) = \frac{\sum\limits_{\kappa \in \mathcal{K}} d_i^{\kappa}(t) \frac{V_i P_i \sigma_i U_i^{\kappa}(t)}{\left(P_i - \sigma_i\right) U_i^{\kappa}(t) + V_i \sigma_i}}{d_i(t)}.
\end{align}
The cell capacity depends on the free flow speeds of each class $U_i^{\kappa}(t)$, as well as on the share of vehicles of each class in the cell, and is lower or equal to the maximum value $Q_i(t) \le Q_i^{\max}$, which is obtained by setting all $U_i^{\kappa}(t) = V_i$, $Q_i^{\max} = V_i \sigma_i$.
Function $F_{i-1}(t)$ models the capacity drop, and $\rho_i(t)=\sum_{\kappa\in\mathcal{K}}\rho_i^{\kappa}(t)$ is the aggregate traffic density.
Where not stated otherwise, the cell parameters will be equal for all cells, and $W = V\frac{\sigma}{P-\sigma}$ yields a triangular fundamental diagram.
The cell length $L$ and time step $T$ are taken so that $L = V T$.  

We prioritize the mainstream flow and only accept on-ramp inflow that the road capacity can support, so a part of vehicles entering the road might have to queue at the on-ramp.
We model the evolution of these queues $n_{i,r}^{\kappa}$, for on-ramps in cell $i$, with
\begingroup
\medmuskip=0.3\medmuskip
\thickmuskip=0.3\thickmuskip
\begin{align}
n_{i,r}^{\kappa}(t) &= n_{i,r}^{\kappa}+\left(\bar{r}_i^{\kappa}(t)-r_i^{\kappa}(t)	\right) T,\\
r_i^{\kappa}(t) &= \begin{cases}
\min\left\{D_{i,r}^{\kappa}(t) , S_{i,r}^{\kappa}(t)\right\},& \kappa \in \mathcal{K}\setminus \mathcal{K}^{*}, \\
\bar{r}_i^{\kappa}(t),& \kappa \in \mathcal{K}^{*},
\end{cases}\\
D_{i,r}^{\kappa}(t) &= \bar{r}_i^{\kappa}(t)+\frac{n_{r,j}^{\kappa}(t-1)}{T},\\
S_{i,r}^{\kappa}(t) &= \frac{n_{i,r}^{\kappa}(t)}{\sum\limits_{m\in \mathcal{K} \setminus \mathcal{K}^{*}}\!\!\!\!\!n_{i,r}^{m}(t)} Q_{i,r}^{\kappa}(t),\\
Q_{i,r}^{\kappa}(t) &= \max\left\{0, \min\left\{S_i(t)-D_i(t),0\right\} - \!\!\!\! \sum\limits_{m \in \mathcal{K}^*} \!\!\! r_i^{m}(t)\right\}.
\end{align}
\endgroup
Here, $\bar{r}_i^{\kappa}(t)$ is the total flow of vehicles arriving at the on-ramp, $S_i(t)$ and $D_i(t)$ are the aggregate supply and demand of cell $i$, and $\mathcal{K}^*$ is the set of prioritized traffic classes.
Vehicles of class $\kappa \in \mathcal{K}^*$ do not queue at the on-ramps, and instead enter the road directly.

One of the benefits of multi-class CTM is that it can exactly define the flow to off-ramps, by representing vehicles with different destinations with different classes.
Let $i$ be a cell with an off-ramp where vehicles of classes $\mathcal{K}_i^r \subset \mathcal{K}$ exit the mainstream.
We may then write
\begin{align}
s_i^{\kappa}(t) &= \begin{cases}
\min\left\{D_i^{\kappa}(t), S_{i+1}^{\kappa}(t), S^{\kappa}_{r,i}(t)\right\},& \kappa \in \mathcal{K}_i^r,\\
0,& \kappa \notin \mathcal{K}_i^r,
\end{cases}\\
S^{\kappa}_{r,i}(t) &= \frac{\rho_i^{\kappa}}{\sum\limits_{m \in \mathcal{K}_i^r} \rho_i^{\kappa}(t)}Q^{\max}_{r,i},
\end{align} 
where $Q^{\max}_{r,i}$ is the capacity of the off-ramp. Finally, we update $q_i^{\kappa}(t)$ accordingly,
\begin{align}
q_i^{\kappa}(t) = \begin{cases}
\min\left\{D_i^{\kappa}(t), S_{i+1}^{\kappa}(t)\right\},& \kappa \notin \mathcal{K}_i^r,\\
0,& \kappa \in \mathcal{K}_i^r.
\end{cases}
\end{align}
 
Out of many ways of modeling capacity drop in first-order traffic models \citep{kontorinaki2017first}, we chose to capture it as a linear reduction of capacity, as in \citep{han2016new}.
Denoting by $\alpha$ the maximum capacity drop ratio under jam traffic density, we have
\begin{equation}
\label{eq:Fcapdrop}
F_i(t) \!=\! W_i\frac{\sigma_{i+1}}{\sigma_i}\left(P_i\!-\!(1\!-\!\alpha)\sigma_i\!-\!\alpha \rho_i(t)\right).
\end{equation}
Note that because of this phenomenon, the actual speed of the congestion wave will be different than $W$.

In our case, vehicles of class $a$ will represent platooned vehicles.
Let there be $\Pi$ platoons, and let platoon $p$ move at speed ${u_p(t) \in \left[U_{\min}, U_{\max} \right]}$, with $U_{\max}<V$.
We denote the position of the platoon head (downstream end) $x_p(t)$, ${x_1(t)> x_2(t) > \dots > x_\Pi(t)}$,  and the reference density of platooned vehicles $\rho_p^*(t)$.
Assuming the length of the platoon is $l_p \ge 2 L$, the traffic density profile in the cells that contain it is
\begingroup
\medmuskip=0.2\medmuskip
\thickmuskip=0.2\thickmuskip
\begin{equation}
\label{eq:MCCTMplat_perfectrho}
\begin{array}{rl}
\rho_{i}^a(t) \!\!\!\!\!\!&= \begin{cases}
0,&\!\!\!\! i<i_\Pi^t(t),\\
\rho_p^*(t)\chi_p^t(t),&\!\!\!\! i = i_p^t(t),\\
\rho_p^*(t),&\!\!\!\! i_p^t(t)<i<i_p^h(t),\\
\rho_p^*(t)\chi_p^h(t),&\!\!\!\! i=i_p^h(t),\\
0,&\!\!\!\! i_p^h(t)<i<i_{p-1}^t(t),
\end{cases}\\
 p \!\!\!\!\!\!&=1, \dots, \Pi,\\
\chi_p^t(t) \!\!\!\!\!\! &= \frac{X_{i_p^t(t)+1}-x_p(t)+l_p}{L},\\
\chi_p^h(t) \!\!\!\!\!\! &= \frac{x_p(t)-X_{i_p^h(t)}}{L},
\end{array}
\end{equation}
\endgroup
where $i_p^h(t) = \left\lceil x_p(t) / L\right\rceil$ and $i_p^t(t) = \left\lceil \left(x_p(t) -l_p\right)/ L\right\rceil$ are the cells in which the platoon head and tail (downstream and upstream end) are, and $i_0^t = N+1$.
The platoon position updates after $T$ will be $x_p(t+1) = x_p(t) + u_p(t) T$, and class $a$ traffic densities need to be updated accordingly, which is achieved by setting
\begingroup
\medmuskip=0.2\medmuskip
\thickmuskip=0.2\thickmuskip
\begin{equation}
\label{eq:Ua_plat}
\begin{array}{rl}
U_{i}^a(t) \!\!\!\!\!\! &= \begin{cases}
V, &\!\! i  <  i_\Pi^t(t),\\
\psi_i^b(t),&\!\! i_p^t(t) \le i < i_p^h(t),\\
\psi_i^h(t), &\!\! i  = i_p^h(t),\\
0, &\!\! i_p^h(t)<i < \frac{i_p^h(t)+i_{p-1}^t(t)}{2}, \\
V, &\!\! \frac{i_p^h(t)+i_{p-1}^t(t)}{2} \le i < i_{p-1}^t(t),
\end{cases}\\
 p\!\!\!\!\!\! &=1, \dots, \Pi,\\
\psi_i^b(t) \!\!\!\!\!\! &= V \frac{ V\rho_p^*(t) - \left(V-U_{i+1}^a(t)\right)\rho_{i+1}^a(t)}{V\rho_i^a(t)},\\
\psi_i^h(t) \!\!\!\!\!\! &= V \left(1 - \left(1-\frac{u_p(t)}{V}\right)\frac{\rho_p^*(t)}{\rho_{i}^a(t)}\right).
\end{array}
\end{equation}
\endgroup
Even if the initial class $a$ density profile differs from the reference, by applying \eqref{eq:Ua_plat} it will converge to \eqref{eq:MCCTMplat_perfectrho}.
Furthermore, the traffic flow overtaking a platoon with density $\rho_p^*$ will be $V \left(\sigma - \rho_p^*\right)$, which is consistent with PDE moving bottleneck models.

Consider a bottleneck at the location of a lane drop, from $n^l_-$ to $n^l_+$ lanes, $n^l_+<n^l_-$.
This corresponds to going from a segment with critical density $\sigma_- = n^l_- \sigma^l$ to $\sigma_+ = n^l_+ \sigma^l$, and the capacity of such bottleneck is $q_b^{\max} =V \sigma_+$.
However, due to the capacity drop phenomenon, in case of excess demand at the bottleneck, its capacity will be decreased once it becomes congested.
A congestion of density $\rho_c$ will be formed, with the density of discharging traffic being $\rho_d$.
The congestion density $\rho_c$ can be calculated from  ${W\left(P_--\rho_c\right) = W\frac{\sigma_+}{\sigma_-}\left( P_- - (1-\alpha)\sigma_- - \alpha\rho_c\right)}$, so that
\begin{align}
\rho_c&= \frac{P_-(\sigma_- - \sigma_+) + (1-\alpha)\sigma_-\sigma_+}{\sigma_- - \alpha\sigma_+}.
\end{align}
We calculate the discharge density from ${V \rho_d = W\left(P_--\rho_c\right)}$:
\begin{align}
\label{eq:rhod}
\rho_d = \frac{\sigma_- \sigma_+ (1-\alpha)}{\sigma_- - \alpha\sigma_+} < \sigma_+.
\end{align}

Since the outflow from the bottleneck is reduced to $q_d=V\rho_d< V \sigma_+$, arriving vehicles will have to wait and their total travel time increase.
We will use the Total Time Spent (TTS), which is the sum of the total time all vehicles spent on the road and the time all vehicles spent queuing at on-ramps, as the performance index of the system,
\begin{align}
\textrm{TTS} = \sum\limits_{t=1}^{t_{\textrm{end}}} \sum\limits_{\kappa \in \mathcal{K}} \sum\limits_{i=1}^{N} \left( \rho_i^\kappa(t) L + n_{i,r}^{\kappa}(t) \right) T.
\end{align}
In order to minimize TTS, we need to keep the demand at the bottleneck as high as possible, while keeping the bottleneck in free flow.

The multi-class CTM introduced in this section can describe fairly complex phenomena, but in order to do this and have a good spatial resolution of the results, we need to use short cells, which makes simulation and prediction less tractable.
For example, if we want to model platoons in traffic, moving at speeds different than $V$, we need the cell length $L$ to be at most half of the platoon length.
Therefore, $L$ will be on the order of magnitude of tens of meters, so we will need a large number of cells to describe any longer highway stretch.
This results in a system with $N \left| \mathcal{K} \right|$ states, where $\left|\mathcal{K}\right|$ is the number of vehicle classes, and up to $N \left| \mathcal{K} \right|$ control inputs if we assume that we can set separate free flow speeds $U_i^\kappa(t)$ for each class and each cell.
In case we want to use this model to predict the outcome of applying some control action, e.g., as a part of optimization-based control, the problem will be intractable due to the large number of states.

However, we may exploit the specific form of the model and the problem to perform state-space reduction without any approximations.
Note that if the considered model is deterministic and $U_i^\kappa(t) = V$ for all classes and cells, assuming the highway was initially in free flow, the only place where we can expect congestion to emerge is at bottlenecks, where $Q^{\max}_{i+1}<Q^{\max}_i$.
Elsewhere, if the road is in free flow, the future traffic density of a cell $\rho_i^{\kappa}(t+j)$ will be equal to the current traffic density of an upstream cell, ${\rho_i^{\kappa}(t+j) = \rho_{i-j}^{\kappa}(t)}$.
Owing to this, we only need to know the initial traffic densities and follow what happens at the bottlenecks, i.e., how the length of their queues evolve in time, to have an accurate view of the full system.
In the following section, we will derive such simplified model, that will then be used to calculate the appropriate control actions and close the loop.

\subsection{Queueing model}
\label{sec:queueing}

\begin{figure*}[t!]
    \centering
    
    \begin{tikzpicture}
    \node[anchor=south west,inner sep=0] (image) at (0,0) {\includegraphics[width=\textwidth]{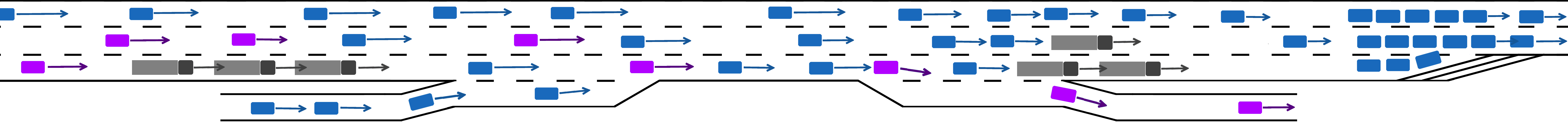}};
    \begin{scope}[x={(image.south east)},y={(image.north west)}]
        \draw[->] (0.16,1.1) -- +(0.08,0) node[above] {\hspace{-1cm}$q^{\rm{out}}_2$};
        \draw[->] (0.67,1.1) -- +(0.08,0) node[above] {\hspace{-1cm}$q^{\rm{out}}_1$};
        \draw[->] (0.92,1.1) -- +(0.08,0) node[above] {\hspace{-1cm}$q^{\rm{out}}_b$};
    \end{scope}
\end{tikzpicture}

    \begin{tikzpicture}[start chain=going right,>=latex,node distance=1cm]

        \node[draw,rectangle,on chain, minimum size=1cm] (n2) {$n_2$};
        \node[draw,rectangle,on chain, minimum size=1cm] (n1) {$n_1$};
        \node[draw,rectangle,on chain, minimum size=1cm] (nb) {$n_b$};

        \draw[->] (n2.east) -- +(1cm,0) node[above] {\hspace{-1cm}$q^{\rm{out}}_2$};
        \draw[->] (n1.east) -- +(1cm,0) node[above] {\hspace{-1cm}$q^{\rm{out}}_1$};
        \draw[->] (nb.east) -- +(1cm,0) node[above] {\hspace{-1cm}$q^{\rm{out}}_b$};
        \draw[<-] (n2.west) -- +(-1cm,0) node[left] {$q^{\rm{in}}_2$};
\end{tikzpicture}

    
    \caption{Queues corresponding to static and moving bottlenecks. The static bottleneck corresponds to $n_b$, the downstream platoon to $n_1$ and the upstream platoon to $n_2$. The overtaking flow of the downstream platoon $q_1^{\rm{out}}$ is limited to one lane of traffic, $q_1^{\rm{cap}} = V {\sigma}_l$, and the overtaking flow of the upstream platoon $q_2^{\rm{out}}$ is limited to two lanes of traffic, $q_2^{\rm{cap}} = V 2 {\sigma}_l$. Both the inflow from the on-ramp and the outflow to the off-ramp will factor in the inflow to the downstream platoon queue $q_1^{\rm{in}}$.}
    \label{fig:queue_road}
\end{figure*}
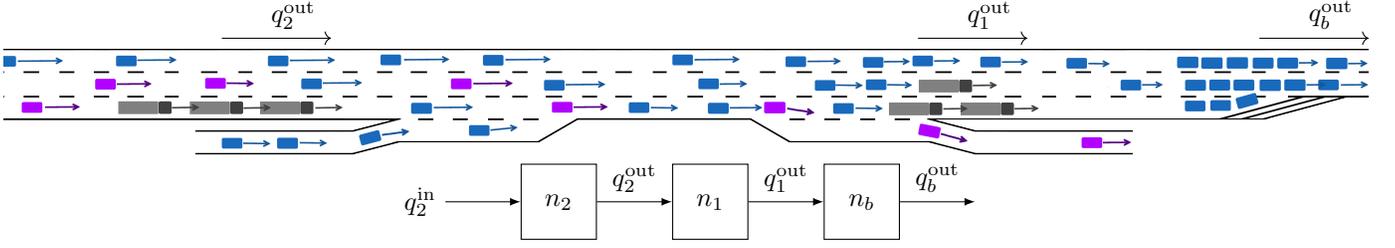

In this work, we study the situation when there is a single bottleneck at the downstream end of the considered stretch of highway, and want to predict its outflow based on the control action we chose for the platoons.
Apart from this stationary bottleneck, platoons themselves can act as moving bottlenecks, since they will be moving slower than the rest of the traffic.
We propose modeling this highway stretch using a queuing-based model, with queue length at the stationary bottleneck $n_b$ and queue lengths at the platoons $n_p$, $p=1, \dots, \Pi$ as the only states.
An example of a traffic situation with its corresponding queuing representation is shown in Figure~\ref{fig:queue_road}.

Since this model is used for predicting the evolution of traffic after some time $t_0$, we assume that the current traffic situation $\rho_i^{\kappa}(t_0)$ is fully known and use this to predict the future values of system states.
We enumerate the platoons that are on the considered highway segment at $t=t_0$, $p=1, \dots, \Pi$, and denote their position at that time $x_p$.
We assume that $t_0 = 0$, and that $t$ represents the prediction time after $t_0$.

The evolution of the queue at the bottleneck is given by
\begin{align}
\label{eq:nb}
\dot{n}_b(t) &= q_b^{\rm{in}}(t) - q_b^{\rm{out}}(t),
\end{align}
where the inflow and the outflow are
\begin{align}
\label{eq:qbin}
q_b^{\rm{in}}(t) &= q_b^u(t) + q_b^V(t),\\
\label{eq:qbout}
q_b^{\rm{out}}(t) &= \begin{cases}
q_b^{\rm{in}}(t),& q_b^{\rm{in}}(t) \le q_b^{\rm{cap}} \land n_b(t) = 0,\\
q_b^{\rm{dis}}, & q_b^{\rm{in}}(t) > q_b^{\rm{cap}} \lor n_b(t)>0.
\end{cases}
\end{align}
Typically, due to capacity drop, the discharge rate of the queue at the bottleneck $q_b^{\rm{dis}}$ will be lower than its capacity $q_b^{\rm{cap}}$, $q_b^{\rm{dis}}<q_b^{\rm{cap}}$.
Mirroring the behaviour of the multi-class CTM, we set ${q_b^{\rm{cap}} = V \sigma_+ = Q^{\max}_+}$ and ${q_b^{\rm{dis}} = V \rho_d}$, according to \eqref{eq:rhod}.
The inflow to the queue at the bottleneck $q_b^{\rm{in}}(t)$ consists of two parts that travel at different speeds.
The first part, $q_b^u(t)$, models the part of the demand that originates from the arrival of the platooned vehicles,
\begingroup
\medmuskip=0.3\medmuskip
\thickmuskip=0.3\thickmuskip
\begin{align}
\label{eq:qbu}
q_b^u(t) &= \begin{cases}
u_p \sigma_l,& t_p^u \le t \le t_p^u+\frac{l_p}{V}, p=1,\dots, \Pi,\\
0,& \text{otherwise},
\end{cases}\\
t_p^u &= \frac{X_b-x_p}{u_p},
\end{align}
\endgroup
where $t_p^u$ represents the time at which platoon $p$ reaches the bottleneck, and the second part consists of the background traffic travelling at free flow speed $V$,
\begingroup
\medmuskip=0.1\medmuskip
\thickmuskip=0.1\thickmuskip
\begin{align}
q_b^V(t) &= \begin{cases}
q_p^{\rm{out}}(\frac{x_p+Vt-X_b}{V-u_p}),&\!\!\!\! \max\left\{t_p^V,t_{p-1}^u\right\}\le t \le t_p^u, \\
V\rho(X_b-Vt),&\!\!\!\! \text{otherwise},
\end{cases}\\
p&=1,\dots,\Pi,\\
t_p^V &= \frac{X_b-x_p}{V}.
\end{align}
\endgroup
Here, the position of the bottleneck is $X_b=x_0$, and $l_p$ is the length of platoon $p$.
We assume that the platoon will approach the bottleneck taking up one lane, thus its density will be equal to the critical density per lane $\sigma_l$.
The second part of the inflow $q_b^V(t)$ originates either from the initial traffic situation,
\begin{align}
\label{eq:rhoinit}
\rho(x) = \!\!\!\sum\limits_{\kappa \in \mathcal{K} \setminus \mathcal{K}^{\Pi}} \!\!\!\rho_i^{\kappa}(0), \quad X_i \le x < X_{i+1},
\end{align}
where $\mathcal{K}\setminus\mathcal{K}^{\Pi}$ is the set of all vehicle classes excluding the platooned vehicles class,
or is the delayed overtaking flow of some platoon.
Each platoon travels at its individual speed $u_p \le V$, and we assume that this speed is constant during the prediction horizon.
Furthermore, we assume that the platoon speeds are such that there is no platoon merging prior to reaching the bottleneck, $t_{p-1}^u > t_p^u$.

Under these assumptions, we define the evolution of the queue at each of the platoons as
\begin{align}
\dot{n}_p(t)& = \frac{V-u_p}{V}\left(q_p^{\rm{in}}(t) - q_p^{\rm{out}}(t)\right) , \quad 0 \le t\le t_p^u,
\end{align}
for $p=1,\dots,\Pi$.
The evolution of queues is defined until time $t_p^u$, when the queue at the platoon is added to the queue at the bottleneck,
\begin{align}
\label{eq:nbtpu}
n_b(t_{p}^u+) = n_b(t_p^u)+n_p(t_p^u).
\end{align}
The outflow and inflow are defined the same way as with the bottleneck queue,
\begin{align}
q_p^{\rm{out}}(t) &= \begin{cases}
q_p^{\rm{in}}(t),& q_p^{\rm{in}}(t) \le q_p^{\rm{cap}}(t) \land n_p(t) = 0,\\
q_p^{\rm{dis}}(t), & q_p^{\rm{in}}(t) > q_p^{\rm{cap}}(t) \lor n_p(t)>0,
\end{cases}\\
q_p^{\rm{in}}(t) &= \begin{cases}
q_{p+1}^{\rm{out}}\left(\frac{(V-u_p)t-x_p+x_{p+1}}{V-u_{p+1}} \right),& t>\frac{x_p-x_{p+1}}{V-u_p},\\
V \rho(x_p-(V-u_p)t,0),& t \le \frac{x_p-x_{p+1}}{V-u_p},
\end{cases}
\end{align}
except here we assume $q_p^{\rm{dis}}(t) = q_p^{\rm{cap}}(t)$, and allow $q_p^{\rm{cap}}(t)$ to vary in time and be used as a control input.
Since the considered road stretch has three lanes, we assume here that platoons can either take one lane or two lanes.
In case platoon $p$ is taking one lane at time $t$, we set $q_p^{\rm{cap}}(t)=V(\sigma_- - \sigma_l)$, and if it is taking two lanes, $q_p^{\rm{cap}}(t) = V(\sigma_- - 2 \sigma_l)$.

The model can be simplified by adopting a coordinate transfer $\tau_p = \frac{x_p-X_b+Vt}{V-u_p}$, $t = \frac{V-u_p}{V}\tau_p+\frac{X_b-x_p}{V}$, for each platoon, which yields
\begin{align}
\frac{\textrm{d} n_p(t(\tau_p))}{\textrm{d} \tau_p} & = q_p^{\rm{in}}(t(\tau_p)) - q_p^{\rm{out}}(t(\tau_p)) , \quad  t_p^V \le \tau_p \le t_p^u
\end{align}
and, taking $\tilde{n}_p(\tau_p) = n_p(t(\tau_p))$, $\tilde{q}_p^{\rm{in}}(\tau_p) = q_p^{\rm{in}}(t(\tau_p))$, and $\tilde{q}_p^{\rm{out}}(\tau_p) = q_p^{\rm{out}}(t(\tau_p))$, we may write
\begin{align}
\label{eq:np}
\dot{\tilde{n}}_p(t) & = \tilde{q}_p^{\rm{in}}(t) - \tilde{q}_p^{\rm{out}}(t) , \quad  t_p^V \le t \le t_p^u
\end{align}
for each $p=1,\dots, \Pi$.
The inflow to the queue at the bottleneck and at platoons can now be simplified to
\begin{align}
q_b^V(t) &= \begin{cases}
\tilde{q}_p^{\rm{out}}(t),& \max\left\{t_p^V,t_{p-1}^u\right\}\le t \le t_p^u,\\
V\rho(X_b-Vt),& \text{otherwise},
\end{cases}\\
\tilde{q}_p^{\rm{in}}(t) &= \begin{cases}
\tilde{q}_{p+1}^{\rm{out}}(t),& t_{p+1}^V<t<t_{p+1}^V,\\
V \rho(X_b-Vt,0),& t \le t_{p+1}^V,
\end{cases}
\end{align}
and the outflow from the platoon becomes
\begin{align}
\label{eq:qpout}
\tilde{q}_p^{\rm{out}}(t) = \begin{cases}
\tilde{q}_p^{\rm{in}}(t),& \tilde{q}_p^{\rm{in}}(t) \le \tilde{q}_p^{\rm{cap}}(t) \land \tilde{n}_p(t) = 0,\\
\tilde{q}_p^{\rm{cap}}(t), & \tilde{q}_p^{\rm{in}}(t) > \tilde{q}_p^{\rm{cap}}(t) \lor \tilde{n}_p(t)>0,
\end{cases}
\end{align}

In case there are on- and off-ramps, their influence can be added to $q_b^V(t)$ and $\tilde{q}_p^{\rm{in}}(t)$.
Denoting $q_k^r(t)$ the inflow from an on-ramp (if ${q_k^r(t)\ge 0}$), or outflow to an off-ramp (if ${q_k^r(t) \le 0}$), we may write 
\begin{align}
\label{eq:qbv}
q_b^V(t) &= q_b^{V\setminus r}(t) + \sum\limits_{k\in K_o^d(t)} \tilde{q}_k^r(t),\\
q_b^{V\setminus r}(t) &=\begin{cases}
\tilde{q}_p^{\rm{out}}(t),& \max\left\{t_p^V,t_{p-1}^u\right\}\le t \le t_p^u,\\
V\rho(X_b-Vt),& \text{otherwise},
\end{cases}\\
K_o^d(t) &= \begin{cases}
K_p^b(t),& \!\!\!\!\max\left\{t_p^V,t_{p-1}^u\right\}\le t \le t_p^u,\\
K_\rho^b(t),& \!\!\!\!\text{otherwise}.
\end{cases}
\end{align}
For the inflow to the queue at platoons, we write
\begin{align}
\label{eq:qpin}
\tilde{q}_p^{\rm{in}}(t) &= \tilde{q}_p^{in\setminus r}(t) + \sum\limits_{k\in K_o^d(t)} \tilde{q}_k^r(t),\\
\tilde{q}_p^{in\setminus r}(t) &= \begin{cases}
\tilde{q}_{p+1}^{\rm{out}}(t),& t>t_{p+1}^V,\\
V \rho(X_b-Vt),& t\le t_{p+1}^V,
\end{cases}\\
K_o^d(t) &= \begin{cases}
K_{p+1}^p(t),& t>t_{p+1}^V,\\
K_\rho^p(t),& t\le t_{p+1}^V.
\end{cases}
\end{align}
Here, $\tilde{q}_k^r(t) = q_k^r(t-\frac{X_b-X_k^r}{V})$, and $K_o^d(t)$ are sets of indices of all on- and off-ramps with positions $X_k^r<X_b$ between the bottleneck or platoon $p$, and the place where their inflows would originate from, 
\begin{align}
K_p^b(t) &= \left\{k \left| x_p^u(t) < X_k^{r} \le X_b, t \ge t_k^r  \right. \right\},\\
K_\rho^b(t) &= \left\{k \left| X_b-Vt < X_k^{r} \le X_b, t \ge t_k^r \right.  \right\},\\
K_{p+1}^p(t) &= \left\{k \left| x_{p+1}^u(t) < X_k^{r} \le x_p^u(t), t \ge t_k^r \right.  \right\},\\
K_\rho^p(t) &= \left\{k \left| X_b-Vt < X_k^{r} \le x_p^u(t), t \ge t_k^r \right.  \right\},\\
t_k^r &= \frac{X_b - X_k^r}{V},
\end{align}
and we define $x_p^u(t)$ as
\begin{align}
x_p^u(t) &= \frac{u_p V t + V x_p - u_p x_p}{V-u_p}.
\end{align}
Note that $q_k^r(t)$ will depend on the local traffic conditions around $X_k^r$ at time $t$.
Furthermore, since a portion of the queue at the platoon will also leave the road via the off-ramp, we reduce $\tilde{n}_p$ at the time when the platoon reaches it,
\begin{align}
\label{eq:npramp}
\tilde{n}_p(t+) = \tilde{n}_p(t) - \Delta_p^{r,k}(t), \quad x_p^u(t) = X_k^r,
\end{align}
and the part of the queue $\tilde{n}_p(t)$ that leaves the highway, $\Delta_p^{r,k}(t)$, depends on the ratio of off-ramp-bound vehicles in the platoon queue.

\begin{figure}[t]
\centering
\includegraphics[width=\columnwidth]{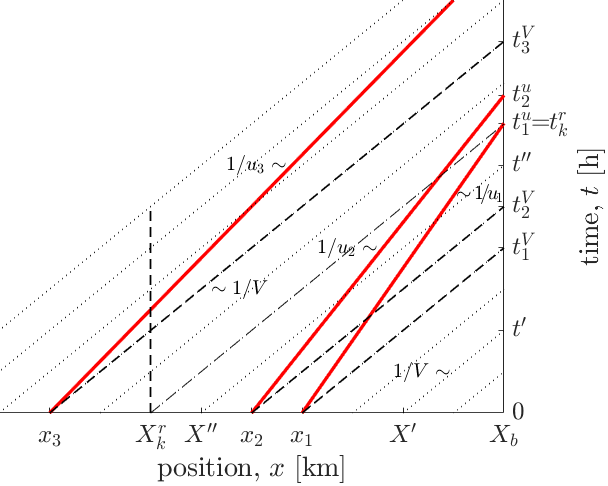}
\caption{Illustration of the queueing model. The dotted lines represent free flow propagation. Platoon trajectories are shown in red. As shown in the figure, at $t = t^{\prime}$, inflow to the bottleneck is ${q_b^{\rm{in}}(t^{\prime}) = V \rho(X^{\prime})}$. At $t = t^{\prime\prime}$, inflow to the bottleneck is ${q_b^{\rm{in}}(t^{\prime\prime}) = \tilde{q}_1^{\rm{out}}(t^{\prime\prime})}$, and inflows to the platoons ${q_1^{\rm{in}}(t^{\prime\prime}) = \tilde{q}_2^{\rm{out}}(t^{\prime\prime})}$, and ${q_2^{\rm{in}}(t^{\prime\prime}) = V \rho(t^{\prime\prime})}$. Ramp $k$ will affect $q_2^{\rm{in}}(t)$ for $\frac{X_b-X_k^r}{V}<t\le t_2^u$, $q_3^{\rm{in}}(t)$ while ${x_3^u(t)\ge X_k^r}$ and ${t<t_3^u}$, and $q_b^{\rm{in}}(t)$ for the rest of time.}
\label{fig:q_illustration}
\end{figure}



An illustration of the derivation of the proposed model is given in Figure~\ref{fig:q_illustration}.
In summary, the proposed model consists of $\Pi+1$ states, whose evolution is described by \eqref{eq:nb} and \eqref{eq:np}.
Inflow to the bottleneck is given by \eqref{eq:qbin}, and consists of the background traffic travelling at free flow speed \eqref{eq:qbv}, and the platoons \eqref{eq:qbu}.
Outflow from the bottleneck is \eqref{eq:qbout}, and there are discontinuous jumps in this state triggered by the arrival of platoons at the bottleneck, \eqref{eq:nbtpu}.
For each platoon queue, inflow is given by \eqref{eq:qpin}, outflow by \eqref{eq:qpout}, and there is a discontinuous jump in the state when the platoon passes an off-ramp, \eqref{eq:npramp}.
The model can be described as a tandem queuing system, with saturation and hysteresis, time-varying structure and jumps.

\section{Control design}
\label{sec:control}

Having defined the simplified model of the system, in this section we will formulate a control law for improving the throughput of the system.
In general, the control objective we consider can be formulated as shaping the traffic flow at some position.
We are looking to maximize the outflow from the bottleneck, which in case there are no off-ramps corresponds to minimizing the total travel time.
In case there are off-ramps the total outflow of the mainstream and of the off-ramps needs to be maximized instead.
We first consider the case when there are no on- or off-ramps and then extend the control to include on- and off-ramps.

The main idea of the proposed control law is to use the platoons as controlled moving bottlenecks whose speed and severity we can control.
We control the moving bottleneck speed by changing the reference speed of the platooned vehicles, and the moving bottleneck severity by changing how many lanes the platoon takes by splitting the platoon and commanding half of the vehicles to drive in parallel to the other half, in the adjacent lane.
By doing this, we are able to first help dissipate the congestion and the static bottleneck, by restricting the flow as much as possible, and then dissipate the congestion in the wake of the moving bottleneck, by reducing the moving bottleneck severity while making sure the static bottleneck remains in free flow.
The proposed control laws rely on the simplified queuing prediction model, and will be described in the remainder of the chapter.


\subsection{Ideal actuation}
\label{sec:ctrl_ideal}

In order to have a baseline for comparing the performance of our proposed control laws, we first consider the ideal case, assuming we can fully control all traffic, and that we can control every class of traffic independently.
This corresponds to having a 100\% penetration rate of connected, communicating, and controlled vehicles, and knowing each vehicle's destination.
Since we already assumed that the demand of off-ramp-bound vehicles is lower than the capacity of the off-ramp, we only need to minimally delay the mainstream-bound background traffic so that the demand at the bottleneck never exceeds its capacity.
This is equivalent to ensuring that the traffic density immediately upstream of the bottleneck $\rho_{i_b}(t) \le \sigma_+$ for all $t$, and can be achieved by setting
\begingroup
\medmuskip=0.1\medmuskip
\thickmuskip=0.1\thickmuskip
\begin{equation}
\label{eq:idealactuate}
\begin{aligned}
\!\!\!\!U^{b}_{i_b}(t) &= V,\\
\!\!\!\!U^{b}_i(t) &= \min\left\{V, \max\left\{ U_{\min}^b,\psi_i^b(t) \right\}\right\}, \quad i = 1, \dots, i_b-1\!\!\!\!\!\!\! \\
\!\!\!\!\psi_i^b(t) &= \frac{V}{\rho_i^b(t)} \left(\rho_i^{b*}(t) - \frac{V -U_{i+1}^b(t)}{V}\rho_{i + 1}^b(t) \right),\\
\!\!\!\!\rho_i^{b*}(t) &= \begin{cases}
\sigma_+ \!\!- \rho_p^*,&\!\!\!\!\!\! \frac{X_b-x_p(t)}{u_p(t)}<\frac{X_b-X_i}{V}<\frac{X_b-x_p(t)+l_p(t)}{u_p(t)}+\frac{L}{V},\!\!\!\!\!\!\! \\
\sigma_+,&\!\!\!\!\!\! \text{otherwise},
\end{cases}\\
\!\!\!\!p &= 1, \dots, \Pi.
\end{aligned}
\end{equation}
\endgroup
This way, the mainstream-bound background traffic is regulated so that the total demand at the bottleneck, including the arriving platoons, is kept as close to its capacity as possible without exceeding it.
The mainstream-bound background traffic is delayed minimally, while the platoons and the off-ramp-bound background traffic experience no delay, travelling at their respective maximum speeds.

\subsection{Platoon-actuation not aware of on- or off-ramps}
\label{sec:ctrl_noramps}


The control objective, maximizing the throughput, i.e., the outflow $q_b^{\rm{out}}$, can be achieved by keeping ${n_b=0}$ and ${q_b^{\rm{in}} = q_b^{\rm{cap}}}$.
Additionally, we require that the queue at the platoon is already discharged when the platoon reaches the bottleneck, $n_p(t_p^u)=0$.
Therefore we employ control law
\begin{align}
\label{eq:qpcap_tilde}
\tilde{q}_p^{\rm{cap}}(t) &= \begin{cases}
q^{\rm{ref}}(t),&  n_b(t)=0 \land t \ge t_{p-1}^u,\\
\tilde{q}_{p-1}^{\rm{cap}}(t), & \tilde{n}_{p-1}(t) = 0 \land t<t_{p-1}^u,\\
Q^{\rm{lo}},& \text{otherwise},
\end{cases}
\end{align}
where the reference flow $q^{\rm{ref}}(t)$ can be externally determined.
For maximizing the throughput, we set
\begin{align}
\label{eq:qref}
q^{\rm{ref}}(t) =Q^{\rm{hi}} - q_b^{u}(t),
\end{align}
taking the largest admissible $Q^{\rm{hi}} \le q_b^{\rm{cap}}$.
In order to compute the current ${q_p^{\rm{cap}}(t) = \tilde{q}_p^{\rm{cap}}(t_p^V)}$ for all platoons, we need to predict $n_b$ until $t_{\Pi}^V$, which requires calculating $q_\Pi^{\rm{cap}}(0)$ and $q_p^{\rm{cap}}(t)$ for ${0\le t \le \min\left\{t_p^u, t_{\Pi}^V\right\}}$.

Assuming this control law is applied, we set the speed of each platoon so that $n_p(t_p^u) = 0$ and ${n_b(t)=0, t_p^c \le t \le t_p^u}$, with minimum $t_p^c$, where
\begin{align}
t_p^c\ge \max\left\{t_p^V, t_{p-1}^{u_{p-1}}+\frac{l_{p-1}}{V}\right\}.
\end{align}
This is achieved when 
\begingroup
\medmuskip=0.3\medmuskip
\thickmuskip=0.3\thickmuskip
\begin{align}
\label{eq:tpu_cond}
\tilde{n}_p(t_p^u) = \tilde{n}_p(t_p^c) + \!\!\int \limits_{t_p^c}^{t_1^u}  \!\! \tilde{q}_1^{\rm{in}}(t) \textrm{d} t - Q^{\rm{hi}}(t_1^u-t_1^c) = 0.
\end{align}
\endgroup
For $p=1$, in case it is known that $t_2^V<t_1^u$, \eqref{eq:tpu_cond} simplifies to
\begin{align}
\tilde{n}_1(t_1^u) &= \tilde{n}_1(t_2^V) + Q^{\rm{lo}} (t_1^u-t_1^c) - Q^{\rm{hi}}(t_1^u-t_1^c) = 0,\\
u_1 &= \frac{\left(Q^{\rm{hi}}-Q^{\rm{lo}}\right)\left(X_b-x_1\right)}{\tilde{n}_1(t_2^V) + \left(Q^{\rm{hi}}- Q^{\rm{lo}} \right) t_1^c},
\end{align}
since we can explicitly calculate 
\[
\tilde{n}_1(t_2^V) = \int_{t_1^V}^{t_2^V} V \rho(X_b-Vt,0)\textrm{d} t-Q^{\rm{lo}}(t_2^V-t_1^V).
\]
Otherwise, $u_p$ is calculated by solving \eqref{eq:tpu_cond} numerically, and can be obtained as a by-product of iterating the prediction steps for $n_b$ and $\tilde{n}_p$.
The simplest way of calculating $u_p$ is to initialize it to
\begin{align*}
\min\left\{U_{\max}, u_{p-1}\frac{X_b-x_p}{X_b-x_{p-1}+l_{p-1}}\right\},
\end{align*}
and then decrease it until either $u_p = U_{\min}$ or \eqref{eq:tpu_cond} is satisfied.
This also ensures that $u_p$ is constrained to be within the range 
\begin{align}
U_{\min} \le u_p \le \min\left\{U_{\max}, u_{p-1}\frac{X_b-x_p}{X_b-x_{p-1}+l_{p-1}}\right\},
\end{align}
which is required for the limitations to be met if there is no platoon merging.

\subsection{Platoon-actuation aware of on- or off-ramps}
\label{sec:ctrl_wramps}

Consider now the case when there are on- or off-ramps.
In order to predict the evolution of queues, which is needed for computing the control inputs, we need to know the ramp flows $\tilde{q}_k^r(t)$ in advance.
This information can be hard to obtain, since it will depend on the routing decisions of individual drivers constituting the background traffic.
Therefore, we use the predicted ramp flows.

If ramp $k$ is an on-ramp, we can replace the actual ramp flow with its average ${\hat{q}_k^r = \bar{q}_k^r}$, which in reality can be determined statistically.
If ramp $k$ is an off-ramp, we can employ the standard assumption that some constant ratio of vehicles $R_k$ leave the road via the off-ramp.
We can then write
\begingroup
\medmuskip=0.3\medmuskip
\thickmuskip=0.3\thickmuskip
\begin{align}
\hat{q}_k^r(t) &= -R_k \left( \tilde{q}^{in,r}_k(t) + \sum\limits_{l\in K_o^{k,r}(t)} \tilde{q}_l^r(t)\right),\\
\tilde{q}^{in,r}_k(t) &= \begin{cases}
 q_b^V(t) ,& x_1^u(t)<x_k^r<X_b\\
 \tilde{q}_{p+1}^{\rm{out}}(t),& x_{p+1}^u(t)<x_k^r<x_p^u(t)
\end{cases}\\
K_o^{k,r}(t) &= \begin{cases}
\left\{ l| x_1^u(t) < x_l^r<x_k^r\right\},& \hspace*{-0.3cm} t>t_1^V, x_k^r < x_{p-1}^u(t)\\
\left\{ l| x_p^u(t) < x_l^r<x_k^r\right\},& \hspace*{-0.3cm}  x_k^r < x_{p-1}^u(t), p>1\\
\left\{ l| X_b-Vt<x_l^r<x_k^r \right\},& \hspace*{-0.3cm}\text{otherwise}
\end{cases}
\end{align}
\endgroup
depending on the origin of the flow to off-ramp $k$ at time $t$.

The portion of queue at platoon $p$ that leaves the highway at off-ramp $k$ can be estimated to be
\begin{align}
\tilde{n}_p(t+) = (1-R_k) \tilde{n}_p(t), \quad  x_p^u(t) = x_k^r,
\end{align}
and we may now apply a control law similar to the one derived for the case when there are no on- and off-ramps.
We modify \eqref{eq:qpcap_tilde} to take into account the fact that there might be some off-ramps $k \in K^*$ whose flow we do not want to obstruct.
Since it is not possible to selectively allow the off-ramp-bound traffic to pass without also releasing the mainstream-bound traffic, we will only allow unrestricted flow towards those off-ramps by setting ${\tilde{q}_p^{\rm{cap}} = Q^{\rm{hi}}}$ if there are other platoons downstream that are regulating the inflow to the bottleneck.
The updated control law is
\begingroup
\medmuskip=0.3\medmuskip
\thickmuskip=0.3\thickmuskip
\begin{align}
\label{eq:qpcap_tilde_onoff}
\tilde{q}_p^{\rm{cap}}(t) &= \begin{cases}
q^{\rm{ref}}(t),& \hspace*{-0.3cm}  n_b(t)=0 \land t \ge t_{p-1}^u,\\
Q^{\rm{hi}},& \hspace*{-0.3cm} K_{p}^{p-1*}(t) \neq \emptyset \land t< t_{p-1}^u,\\
\tilde{q}_{p-1}^{\rm{cap}}(t), & \hspace*{-0.3cm} K_{p}^{p-1*}(t) = \emptyset \land \tilde{n}_{p-1}(t) = 0 \land t<t_{p-1}^u,\\
Q^{\rm{lo}},& \hspace*{-0.3cm} \text{otherwise},
\end{cases}
\end{align}
\endgroup
where $K_{p}^{p*}(t) = K_{p+1}^{p}(t) \cap K^*$.

The platoon speeds are again obtained in the course of predicting the queue evolution, as described in the previous subsection.


\section{Analysis}
\label{sec:analysis}

In order to understand the effects and limitations this control law will have in realistic situations, we first study it under simplified conditions, in an idealised situation.
Whereas in simulations the inflow of background traffic will vary in time and take random values belonging to some range, and platoons arrive with exponentially distributed gaps, we first assume constant background traffic inflow ${Q^{\rm{in}}(t) = Q^{\rm{in}}}$ and periodic platoon arrivals, with period  $\tau_{\pi}$, and each platoon consisting of $n_{\pi}$ passenger car equivalents, and then allow the inflow and gaps between two platoons vary within some range.
In this section, we derive:
\begin{enumerate}
\item Exact limits on the maximum initial excess congestion for which the controlled system is stable, assuming constant inflow and periodic platoon arrivals,
\item The number of controlled platoons required to fully dissipate the congestion at a static bottleneck and return the road to the unperturbed free flow state, and
\item An estimate of throughput given varying inflow and gap between platoons, i.e., the average inflow for which we are able to dissipate the congestion at the bottleneck with a predefined probability.
\end{enumerate}

The bottleneck will be considered to have capacity $q_b^{\rm{cap}}$, which is reduced to $q_b^{\rm{dis}}$ in case there is capacity drop, ${q_b^{\rm{dis}} < q_b^{\rm{cap}}}$.
We study the case when the bottleneck is already congested at initial time.
If the platoon arrives at a congested bottleneck, its vehicles are added to the bottleneck queue.
Otherwise, if there is no queue at the platoon and it arrives at a bottleneck in free flow, the platoon passes through the bottleneck without causing traffic breakdown.

In summary, the system that we study in this section is
\begin{equation}
\label{eq:modelsummary}
\begin{aligned}
\dot{n}_b(t) \! &= q_b^{\rm{in}}(t) - q_b^{\rm{out}}(t),\\
q_b^{\rm{in}}(t) \! &= \!\begin{cases}
\tilde{q}_p^{\rm{out}},& \max \left\{ t_p^V, t_{p-1}^u\right\}\le t \le t_p^u,\\
Q^{\rm{in}}(t),& \textrm{otherwise},
\end{cases}\\
q_b^{\rm{out}}(t) \! &= \!\begin{cases}
q_b^{\rm{in}}(t),& q_b^{\rm{in}}(t) \le q_b^{\rm{cap}} \land n_b(t) = 0,\\
q_b^{\rm{dis}},& q_b^{\rm{in}}(t)> q_b^{\rm{cap}} \lor n_b(t)>0,
\end{cases}\\
\dot{\tilde{n}}_p(t) \! &= \tilde{q}_p^{\rm{in}}(t) - \tilde{q}_p^{\rm{out}}(t),\\
\tilde{q}_p^{\rm{in}}(t) \! &= \!\begin{cases}
\tilde{q}_{p+1}^{\rm{out}},& t_{p+1}^V<t<t_{p+1}^u,\\
Q^{\rm{in}}(t),& t \le t_{p+1}^V,
\end{cases}\\
\tilde{q}_p^{\rm{out}}(t) \! &= \!\begin{cases}
\tilde{q}_p^{\rm{in}}(t),& \tilde{q}_p^{\rm{in}}(t) \le \tilde{q}_p^{\rm{cap}}(t) \land \tilde{n}_p(t) = 0,\\
\tilde{q}_p^{\rm{dis}},& \tilde{q}_p^{\rm{in}}(t) > \tilde{q}_p^{\rm{cap}}(t) \lor \tilde{n}_p(t)>0,
\end{cases}\\
n_b(t_p^u\!+\!) \! &= \!\begin{cases}
\!n_b(t_p^u)\!+\!n_p(t_p^u)\!+\!n_{\pi},& \!\!\! n_b(t_p^u)\! +\! \tilde{n}_p(t_p^u)\!>\!0,\!\!\!\!\!\!\\
\!0,& \!\!\! n_b(t_p^u)\!+\!\tilde{n}_p(t_p^u)\!=\!0,\!\!\!\!\!\!
\end{cases}\\
p \! &= 1, \dots, \Pi
\end{aligned}
\end{equation}
where $\tilde{q}_p^{\rm{cap}}(t)$ is governed by control law \eqref{eq:qpcap_tilde}.

\subsection{Constant inflow and periodic platoon arrivals}

We study the stability of the queue at the bottleneck under conditions of constant inflow and periodic platoon arrivals for different initial bottleneck queue lengths.
First, in case no control is applied, the system is stable if
\begin{equation}
\label{eq:uncont_stable}
\tilde{Q}^{\rm{in}} = Q^{\rm{in}} + \frac{n_{\pi}}{\tau_{\pi}} < q_b^{\rm{dis}},
\end{equation}
i.e., if the average total inflow is less than the dissipating flow of the bottleneck, and its queue length will go to zero regardless of its initial value.

However, if the platoons can be controlled, we are able to extend the range of $Q^{\rm{in}}$ for which the system is stable.
In this case, it is of interest to study what is the maximum initial queue length $\bar{n}_b^0$ for which the system is stable for every given $Q^{\rm{in}}$.
We control the platoons by setting their speeds and how many lanes they take.
The platoon speed can range from some set minimum speed $U^{\min}$ to the free flow speed of all traffic $V$.
By setting the number of lanes a platoon occupies, we control the maximum overtaking flow at its position, with overtaking flow of $Q^{\rm{lo}}$ corresponding to maximum number of lanes taken, and $Q^{\rm{hi}}$ corresponding to one lane taken.

The case we are considering assumes that the flow values are arranged as
\begin{equation}
\label{eq:flow_ordering}
Q^{\rm{lo}} < q_b^{\rm{dis}} < Q^{\rm{in}} < \tilde{Q}^{\rm{in}} < Q^{\rm{hi}} \le q_b^{\rm{cap}},
\end{equation}
and the uncontrolled system is unstable.
The length of the considered road segment is $\ell$ and a platoon moving at speed $u_k$ traverses it and reaches the bottleneck after
\[
\tau_p^k = \frac{\ell}{u_k}.
\]
We can define the initial queue length ${\bar{n}_b^0 = n_b(\ell/V)}$ as the queue length at the bottleneck at the time when the overtaking flow from the platoon reaches it.
This is equivalent to saying that there is $\bar{n}_b^0$ excess congestion to be dissipated, signifying how many vehicles need to be temporarily removed from the inflow in order for the bottleneck to return to free flow.
For the first platoon entering the road segment, the entire congestion will be located at the bottleneck, and for subsequent platoons, the initial excess congestion $\bar{n}_b^k$ will be distributed between the bottleneck and downstream platoons that previously entered the road. 
The system is stable if ${\bar{n}_b^{k+1}<\bar{n}_b^k}$, i.e., if every subsequent platoon has less excess congestion to dissipate.
Since we are looking for maximum $\bar{n}_b^k$ for which this holds, we study the situation when maximum control action is applied, corresponding to setting the platoon speed to $u_k = U^{\min}$ with maximum overtaking flow $Q^{\rm{lo}}$ until the queue at the bottleneck is dissipated, which happens after
\[
\tau_c^k = \frac{\bar{n}_b^k}{q_b^{\rm{dis}} - Q^{\rm{lo}}}.
\]
Moving at minimum speed, a platoon will reach the bottleneck after
\[
\tau_p^{\max} = \frac{\ell}{U^{\min}},
\]
so a necessary condition to be able to begin dissipating the congestion is that ${\tau_p^{\max}>\tau_c^k}$, which yields
\begin{equation}
\bar{n}_b^0<\frac{q_b^{\rm{dis}}-Q^{\rm{lo}}}{U^{\min}}\ell.
\end{equation}

%


%


The process of dissipating excess congestion can be split into two phases.
In the first phase, maximum control action is applied and there is a queue at the platoons when they reach the bottleneck.
Given excess congestion $\bar{n}_b^k$, the excess congestion left for platoon $k+1$ to dissipate can be written as
\begin{equation}
\label{eq:barn_update}
\bar{n}_b^{k+1}=a\bar{n}_b^k + b,
\end{equation}
where
\begin{align}
\label{eq:a}
a &= \frac{Q^{\rm{hi}} - Q^{\rm{lo}}}{q_b^{\rm{dis}} - Q^{\rm{lo}}} > 1, \\
\label{eq:b}
b &= \tau_{\pi}\left(Q^{\rm{in}} - q_b^{\rm{dis}}\right) + n_{\pi} - \tau_p^{\max} \left( Q^{\rm{hi}} - q_b^{\rm{dis}} \right) < 0.
\end{align}
Therefore, the stability condition of \eqref{eq:barn_update} is
\[
\bar{n}_b^{k} < \frac{b}{1-a}.
\]
Starting with $\bar{n}_b^0$, we can calculate $\bar{n}_b^k$ by recursing \eqref{eq:barn_update},
\[
\bar{n}_b^k = a^k \bar{n}_b^0 + \sum\limits_{i=0}^{k-1}a^i b.
\]

The second phase begins with platoon $k$ able to dissipate all excess congestion and reach the bottleneck with no queue.
This happens when $\bar{n}_b^k < c$, where 
\begin{equation}
\label{eq:second_phase_cond}
c = \left(q_b^{\rm{dis}} - Q^{\rm{lo}}\right) \left( \tau_p^{\max} - \tau_{\pi}\frac{Q^{\rm{in}} - Q^{\rm{lo}}}{Q^{\rm{hi}}-Q^{\rm{lo}}} \right),
\end{equation}
and given $\bar{n}_b^0$, the transition into the second phase of congestion dissipation after $k'$ platoons, where $k'$ is the lowest integer such that
\[
a^{k'} \bar{n}_b^0 + \sum\limits_{i=1}^{k'-1}a^i b \le c.
\]
Since in this phase the platoons will be reaching the bottleneck that is in free flow, and with no queues, in this phase we may approximate $n_{\pi} \approx 0 $. 

The minimum time when platoon $k$ can reach the bottleneck with no queue is
\begin{equation}
\label{eq:taupk}
\tau_p^k = \frac{\bar{n}_b^k}{q_b^{\rm{dis}} - Q^{\rm{lo}}} + \tau_{\pi} \frac{Q^{\rm{in}} - Q^{\rm{lo}}}{Q^{\rm{hi}} - Q^{\rm{lo}}},
\end{equation}
in which case the platoon is travelling at speed
\[
u_k = \frac{l}{\tau_{\pi} \frac{Q^{\rm{in}} - Q^{\rm{lo}}}{Q^{\rm{hi}} - Q^{\rm{lo}}}+ \frac{\bar{n}_b^k}{q_b^{\rm{dis}} - Q^{\rm{lo}}}}.
\]
Substituting \eqref{eq:barn_update} into \eqref{eq:taupk} we get
\[
\tau_p^{k+1} = \tau_p^{k} - \tau_{\pi} \frac{Q^{\rm{hi}} - Q^{\rm{in}}}{Q^{\rm{hi}} - Q^{\rm{lo}}},
\]
and the speed of platoon $k+1$ will be
\[
u_p^{k+1} = \frac{l}{\frac{l}{u_k} - \tau_{\pi} \frac{Q^{\rm{hi}} - Q^{\rm{in}}}{Q^{\rm{hi}} - Q^{\rm{lo}}}}.
\]
%
The traffic will return to the unperturbed state after platoon $k''$ for which
\[
\tau_p^{k''} = \frac{l}{u_p^{k''}} \le \tau_{\pi}, 
\]
or equivalently,
\[
\bar{n}_b^{k'} \le \left(q_b^{\rm{dis}} - Q^{\rm{lo}}\right) \left( \tau_p^{\min} - \tau_{\pi}\frac{Q^{\rm{in}} - Q^{\rm{lo}}}{Q^{\rm{hi}}-Q^{\rm{lo}}} \right).
\]
Given $\tau_p^{k'}<\tau_p^{\max}$, we can calculate $k''$ by rounding up
\[
k'' = \left\lceil \frac{\tau_p^{k'} - \tau_{\pi}}{\tau_{\pi} \frac{Q^{\rm{hi}} - Q^{\rm{in}}}{Q^{\rm{hi}} - Q^{\rm{lo}}}}\right\rceil.
\]

We summarize the analysis in this proposition:

\begin{prop}
\label{prop1}
Assuming constant inflow $Q^{\rm{in}}(t) = Q^{\rm{in}}$, periodic arrival of platoons with period $\tau_{\pi}$ and ordering of flow values \eqref{eq:flow_ordering}, the queue length $n_b(t)$ of system \eqref{eq:modelsummary} controlled by control law \eqref{eq:qpcap_tilde} is stable and will remain $0$ after some time $t$, if the initial queue length satisfies
\begin{equation}
\bar{n}_b^0<\frac{b}{1-a},
\end{equation} 
where $a$ and $b$ are given by \eqref{eq:a} and \eqref{eq:b}, respectively, and $\bar{n}_b^0 = n_b(\frac{\ell}{V})$.
Furthermore, if this condition is satisfied , the system returns to the unperturbed state with $n_b(t)=0$ and $\tilde{n}_p(t)=0$ after platoon $k''$ reaches the bottleneck.
\end{prop}

\subsection{Varying inflow and platoon arrivals}

The uncertainty coming from the varying inflow of background traffic and random platoon arrivals can be modelled by adding another term to \eqref{eq:barn_update}:
\begin{equation}
\label{eq:barn_update_delta}
\bar{n}_b^{k+1} =a\bar{n}_b^{k} + b + \delta^{k},
\end{equation}
where $\delta^{k} = \delta_{\tau_\pi}^{k} (Q^{\rm{in}}-q_b^{\rm{dis}})+\tau_{\pi}\delta_{Q^{\rm{in}}}^k+\delta_{\tau_{\pi}}^k\delta_{Q^{\rm{in}}}^k$, $\delta_{\tau_\pi}^{k}$ is the difference of the gap between platoon $k-1$ and $k$ from $\tau_{\pi}$, and $\delta_{Q^{\rm{in}}}^k$ is the difference of the average inflow from $Q^{\rm{in}}$ during that time.
We may also write
\[
\bar{n}_b^{k} = a^k \bar{n}_b^0 + \sum\limits_{i=0}^{k-1}a^{k-1-i}\left(\delta^i+ b\right).
\]
\begin{prop}
Assuming that ${\left|\delta^k\right|<\Delta<\left|b\right|}$, if for any $k$ we have
\begin{equation}
\label{eq:nk_stable}
\bar{n}_b^k < \frac{b+\Delta}{1-a},
\end{equation}
with $a$ given by \eqref{eq:a} and $b$ by \eqref{eq:b}, then system \eqref{eq:barn_update_delta} is stable, so we are able to dissipate the congestion at the bottleneck.
Conversely, if for any $k$ we have
\begin{equation}
\label{eq:nk_unstable}
\bar{n}_b^k > \frac{b-\Delta}{1-a},
\end{equation}
then system \eqref{eq:barn_update_delta} is unstable, so in that case the queue at the bottleneck will grow unbounded.
\end{prop}
Consequently the conclusions about stability can be extended to system \eqref{eq:modelsummary} if a suitable bound on uncertainty $\Delta$ can be derived.

For the initial excess congestion between these two values, ${\frac{b+\Delta}{1-a}< \bar{n}_b^0< \frac{b-\Delta}{1-a}}$, $\bar{n}_b^k$ will almost surely satisfy either condition \eqref{eq:nk_stable} or \eqref{eq:nk_unstable} for some $k$, after which the queue stability does not depend on $\delta^k$.
Assuming uniformly distributed $\delta^k$, with $\mathrm{E}\left\{\delta^k\right\}=0$, $\mathrm{Var}\left\{\delta^k\right\} = \frac{\Delta^2}{3}$ (e.g., if $\delta^k \sim \mathcal{U}[-\Delta, \Delta]$), the probability of $\bar{n}_b^k$ satisfying \eqref{eq:nk_unstable} (i.e., failing to decongest the bottleneck) closely follows the logistic curve depending on $\bar{n}_b^0$,
\begin{equation}
\label{eq:punstab}
\mathcal{P}_{uns}(\bar{n}_b^0) \approx \left(1+\exp\left(\frac{\frac{b}{1-a}-\bar{n}_b^0}{\frac{\Delta}{4}}\right)\right)^{-1}
\end{equation}
and the probability of $\bar{n}_b^k$ satisfying \eqref{eq:nk_stable} for some $k$ (i.e., successfully decongesting the bottleneck) is ${\mathcal{P}_{sta}(\bar{n}_b^0) = 1-\mathcal{P}_{uns}(\bar{n}_b^0)}$.

Finally, we may define the estimate of throughput of the controlled system as the maximum $Q^{\rm{in}}$ for which the control algorithm is able to decongest the bottleneck with probability $\mathcal{P}^*$, given an appropriately chosen $\bar{n}_b^0$.
This yields the bound on $Q^{\rm{in}}$,
\begingroup
\medmuskip=0.1\medmuskip
\thickmuskip=0.1\thickmuskip
\begin{equation}
\label{eq:ctrl_stab0}
\tilde{Q}^{\rm{in}}_{\bar{n}_b^0, \mathcal{P}^*\!\!, \Delta} = q_b^{\rm{dis}}\!\!+\frac{Q^{\rm{hi}}\!-q_b^{\rm{dis}}}{\tau_{\pi}}\!\left( \!\!\tau_p^{\max} \!\!- \frac{\bar{n}_b^0 + \frac{\Delta}{4} \log\!\left(\frac{\mathcal{P}^*\!\!}{1-\mathcal{P}^*\!\!}\right) }{q_b^{\rm{dis}}\!-Q^{\rm{lo}}} \!\!\right).
\end{equation}
\endgroup
One suitable choice for $\bar{n}_b^0$ is to take
\begingroup
\medmuskip=0.1\medmuskip
\thickmuskip=0.1\thickmuskip
\begin{equation}
\label{eq:nb0_c}
\scalebox{0.95}{$
\bar{n}_b^0 = \left(q_b^{\rm{dis}}\!\!-Q^{\rm{lo}}\right)\!\left(\tau_p^{\max}\!\!-\tau_{\pi}\right)+n_{\pi}+\frac{Q^{\rm{hi}}\!\!-q_b^{\rm{dis}}}{q_b^{\rm{dis}}-Q^{\rm{lo}}}\frac{\Delta}{4}\!\log\!\left(\!\frac{\mathcal{P}^*\!}{1-\mathcal{P}^*\!}\!\right),
$}
\end{equation}
\endgroup
in which case we have
\begingroup
\medmuskip=0.1\medmuskip
\thickmuskip=0.1\thickmuskip
\begin{equation}
\label{eq:ctrl_stab}
\scalebox{0.92}{$
\!\!\!\tilde{Q}^{\rm{in}}_{\mathcal{P}^*\!\!, \Delta} = Q^{\rm{hi}} - \frac{Q^{\rm{hi}}-q_b^{\rm{dis}}}{q_b^{\rm{dis}}-Q^{\rm{lo}}}\left(\frac{n_{\pi}}{\tau_{\pi}} + \frac{Q^{\rm{hi}}-Q^{\rm{lo}}}{q_b^{\rm{dis}}-Q^{\rm{lo}}} \frac{\Delta}{4} \frac{\log\left(\frac{\mathcal{P}^*\!\!}{1-\mathcal{P}^*\!\!}\right)}{\tau_{\pi}}\right),
$}
\end{equation}
\endgroup
and thus $\bar{n}_b^0$ is equal to $c$ from condition \eqref{eq:second_phase_cond} with ${Q^{\rm{in}} = \tilde{Q}^{\rm{in}}_{\mathcal{P}^*\!\!, \Delta}-\frac{n_{\pi}}{\tau_{\pi}}}$.

\begin{figure*}[t!]
\centering
\begin{subfigure}[t]{0.49\linewidth}
\includegraphics[width=\textwidth]{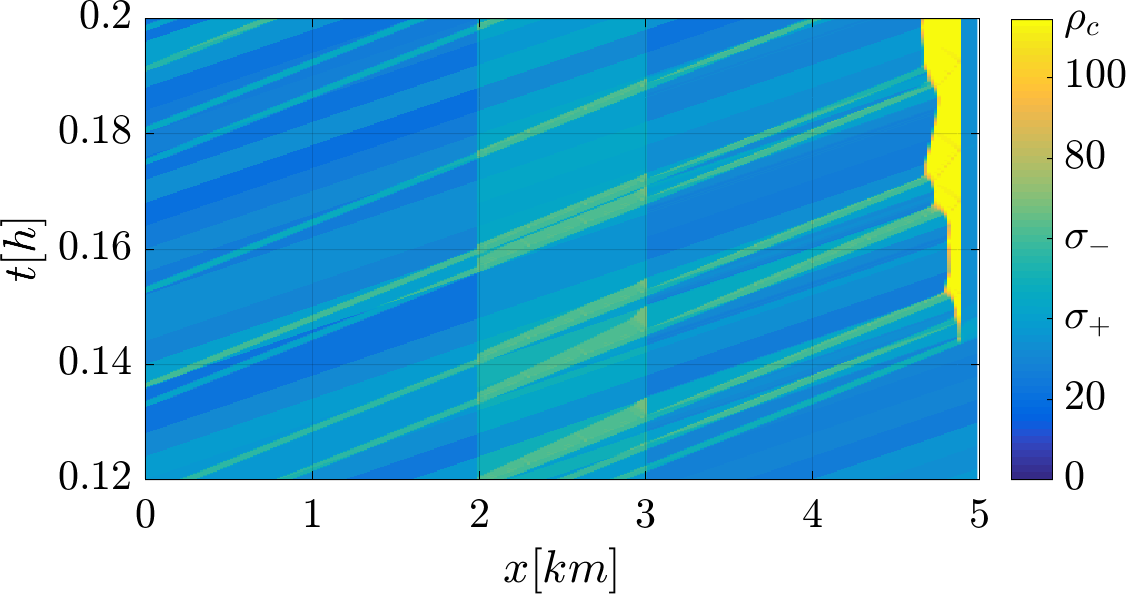}
\caption{No control}
\label{fig:ex_noctrl}
\end{subfigure}
\begin{subfigure}[t]{.49\linewidth}
\includegraphics[width=\textwidth]{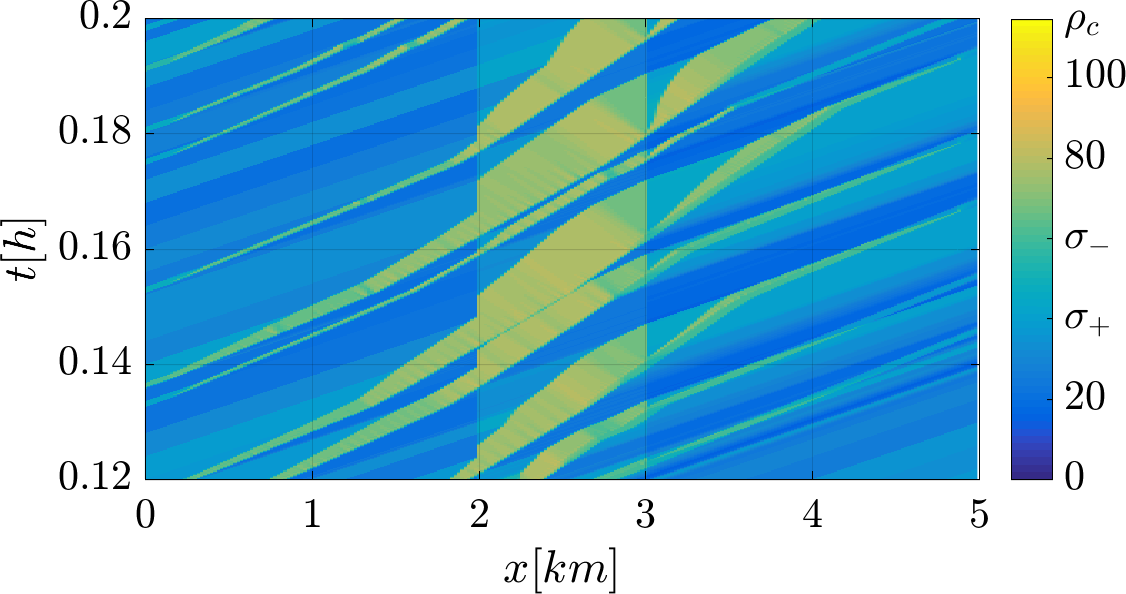}
    \caption{Platoon-actuated control ignoring on- and off-ramps}
    \label{fig:ex_ctrlnoramp}
\end{subfigure}
\\
\noindent
\begin{subfigure}[t]{.49\linewidth}
\includegraphics[width=\textwidth]{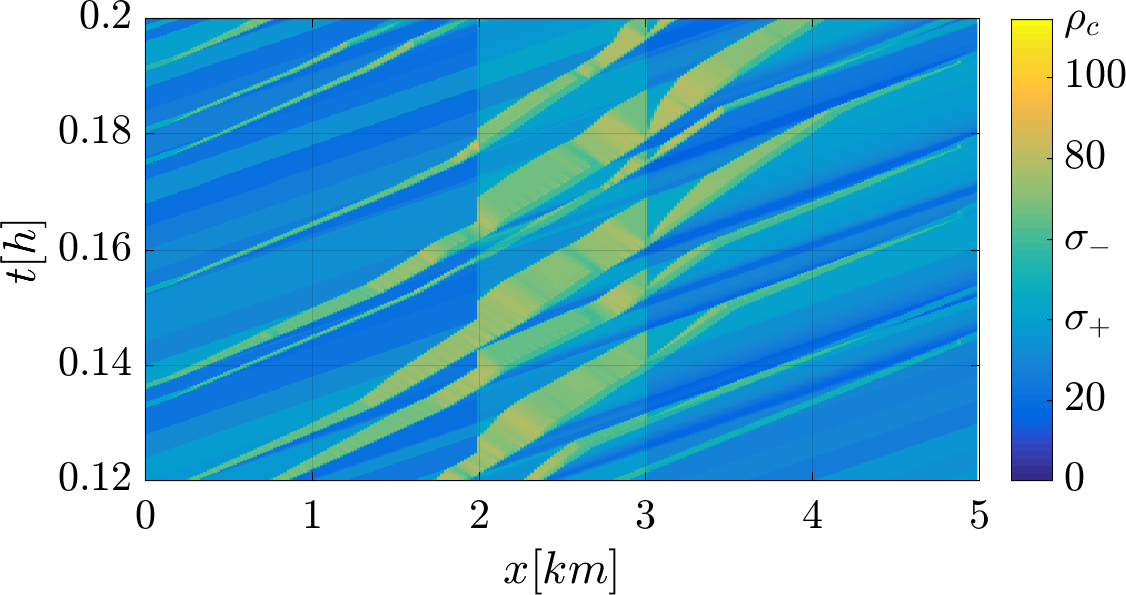}
    \caption{Platoon-actuated control taking on- and off-ramps into account}
    \label{fig:ex_ctrlwramp}
\end{subfigure}
\begin{subfigure}[t]{.49\linewidth}
\includegraphics[width=\textwidth]{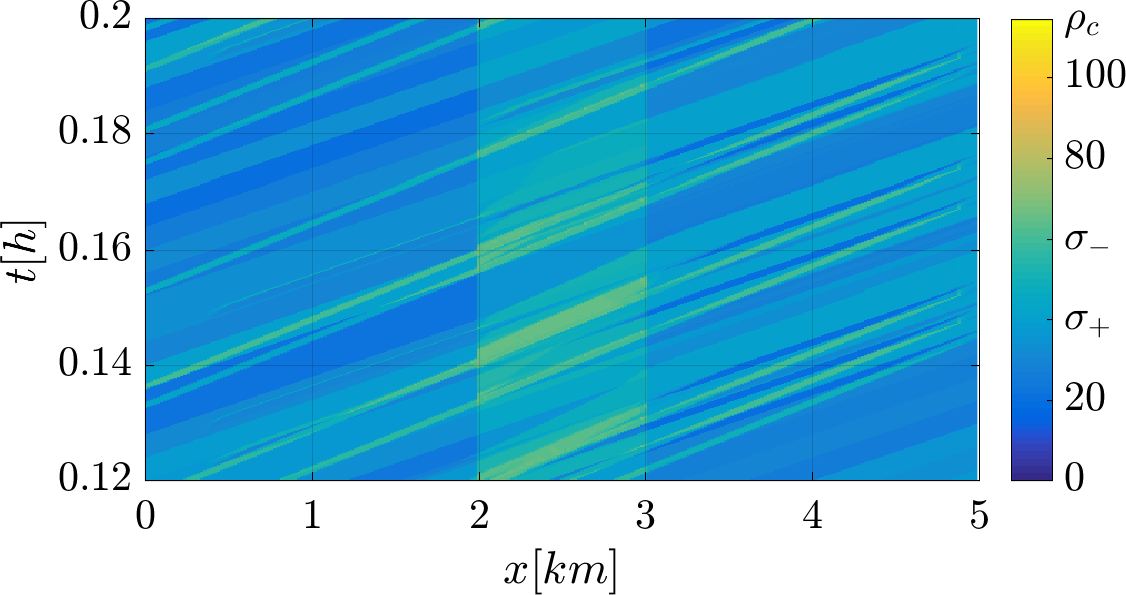}
    \caption{Ideally actuated control}
    \label{fig:ex_ctrlideal}
\end{subfigure}
\caption[]{An example comparing the outcome of the four simulation cases. Traffic density is color-coded, with warmer color representing higher density.}
\label{fig:example}
\end{figure*}

\section{Simulation-based Validation}
\label{sec:sim}
In order to assess the performance of the proposed control laws, we conducted a number of simulation runs, results of which will be presented in this section.
We demonstrate the control laws' effectiveness and show that we are able to eliminate $52.7\%$ of the total delay due to congestion experienced by all vehicles on average, or $75.6\%$ of the total delay by median, compared to the case when no control is applied.


The simulations were executed on a $5$ km long stretch of highway, illustrated by Figure~\ref{fig:queue_road}, with an on-ramp around the $2$ km mark, and an off-ramp around the $3$ km mark.
Most of the highway stretch has three lanes, corresponding to a critical density  of $\sigma_- = 60$ veh/km and capacity of $Q^{\max}_- = 6000$ veh/h, with free flow speed of $V = 100$ km/h.
There is a bottleneck caused by a lane drop $80$ m upstream of the end of the considered stretch, with capacity of $Q^{\max}_+ = 4000$ veh/h.
The capacity drop phenomenon is modelled with $\alpha = 0.4$, which causes the bottleneck capacity to be reduced to $Q^{\rm{dis}}_+ = 3273$ veh/h, representing a $18.2\%$ capacity drop for this road configuration.

We considered three classes of traffic: class $a$ consists of the platoons we control, class $b$ is the mainstream-bound background traffic, and class $c$ the off-ramp-bound background traffic.
The arrival of class $a$ vehicles is modelled as Poisson process with Poisson arrival rate of $\lambda = 81$ platoon/h, $\tau_{\pi}=0.0123$ h. 
We assume that each platoon consists of $2$ passenger car equivalents, although in reality, due to having shorter inter-vehicular gaps, these platoons might be up to about five passenger cars or about three trucks long.
This effect was not included in calculating the TTS, and including it would only further emphasize the benefits of proposed control. 
The inflow of background traffic is assumed to be time-varying and uniformly distributed, changing every $14.4$ seconds.
At the beginning of the highway segment, the demand of mainstream-bound background traffic takes values in $\bar{r}_1^{b}(t)\sim \mathcal{U}(1000, 2000)$ veh/h, and the demand of off-ramp bound traffic is $\bar{r}_1^{c}(t)\sim \mathcal{U}(750, 1250)$ veh/h.
Since the on-ramp and off-ramp are reasonably close, we assume that none of the vehicles entering the highway via the on-ramp will exit it via the off-ramp, $\bar{r}_{i_{\rm{on}}}^{c}(t)=0$ veh/h.
The demand of mainstream-bound traffic at the on-ramp is modelled as $\bar{r}_{i_{\rm{on}}}^{b}(t)\sim \mathcal{U}(900, 1500)$ veh/h.

With the parameters specified in previous paragraph, we may calculate an estimate of the throughput that we may achieve by applying the presented control law.
Using \eqref{eq:ctrl_stab} with $\Delta \!=\! \tau_{\pi}\!\!\left(\max\left(\bar{r}_1^b+\bar{r}_{i_{\rm{on}}}^b\right)\!-\!\mathrm{E}\!\!\left\{\bar{r}_1^b+\bar{r}_{i_{\rm{on}}}^b\right\}\right)$ and $\mathcal{P}^* \!\!=\! 0.9$, we estimate that the throughput would be improved from $\tilde{Q}^{\rm{in}}_{\rm{unc}} = 3273$ veh/h to $\tilde{Q}^{\rm{in}}_{\mathcal{P}^*, \Delta} = 3513.2$ veh/h.
Note that in deriving \eqref{eq:ctrl_stab} we do not take into account the existence of the on-ramp.

The duration of each simulation run is $2$ hours, of which the background traffic inflow is halved for the first $3$ minutes, in order to properly initialize the system, and for the last $12$ minutes, in order to allow the traffic to return to free flow and ensure fair comparison between different control laws.
Simulations are done with four cases of control:
\begin{enumerate}[label=(\alph*)]
\item No control,
\item Platoon-actuated control ignoring on- and off-ramps from Subsection~\ref{sec:ctrl_noramps}, with $\tilde{q}_p^{\rm{cap}}(t)$ given by \eqref{eq:qpcap_tilde},
\item Platoon-actuated control taking on- and off-ramps into account from Subsection~\ref{sec:ctrl_wramps}, with $\tilde{q}_p^{\rm{cap}}(t)$ given by \eqref{eq:qpcap_tilde_onoff}, and
\item Ideally actuated control from Subsection~\ref{sec:ctrl_ideal}, with $U_i^b(t)$ given by \eqref{eq:idealactuate}.
\end{enumerate}
In order to demonstrate the effect applying these control laws has on the traffic, a part of one simulation run is shown in Figure~\ref{fig:example}.




\renewcommand{\thesubfigure}{\alph{subfigure}}
\begin{figure*}[t!]
\centering
\begin{subfigure}[t]{.245\linewidth}
\includegraphics[height=4cm]{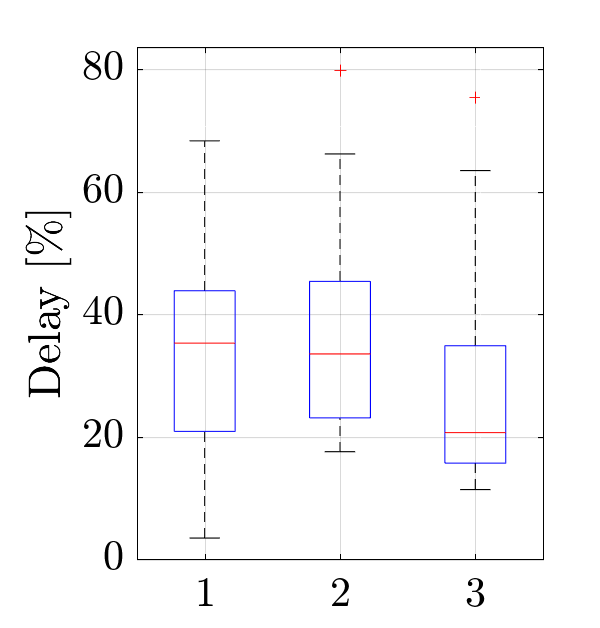}
\caption{Class $a$}
\label{fig:box_noctrl}
\end{subfigure}
\begin{subfigure}[t]{.245\linewidth}
\includegraphics[height=4cm]{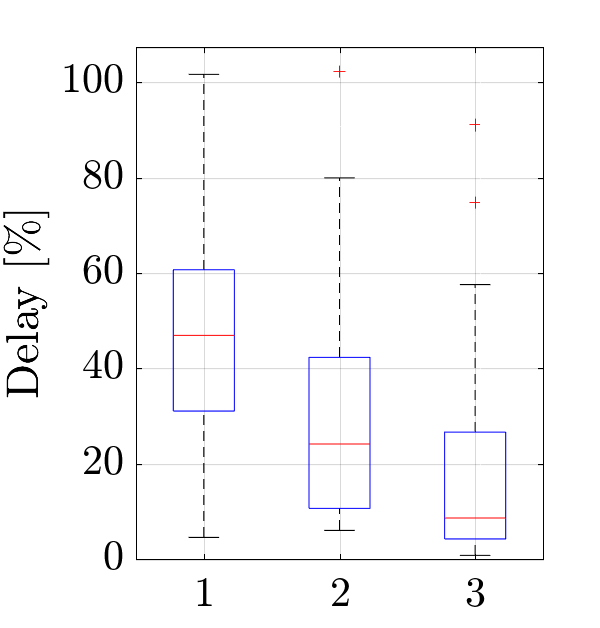}
    \caption{Class $b$}
    \label{fig:box_ctrlnoramp}
\end{subfigure}
\begin{subfigure}[t]{.245\linewidth}
\includegraphics[height=4cm]{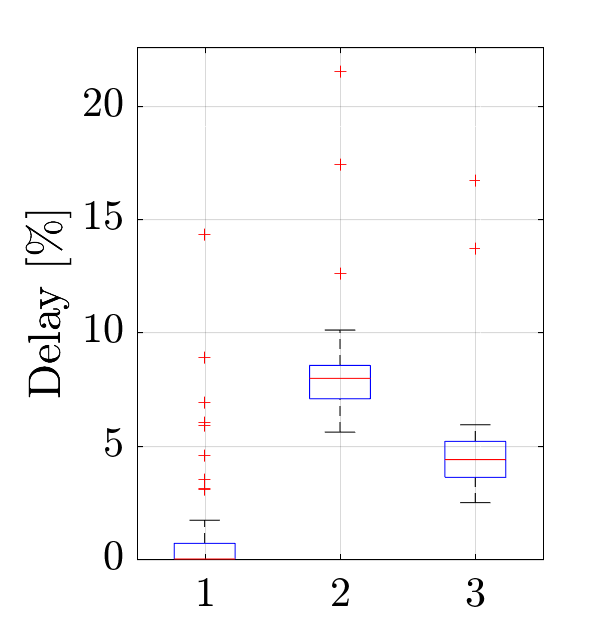}
    \caption{Class $c$}
    \label{fig:box_ctrlwramp}
\end{subfigure}
\begin{subfigure}[t]{.245\linewidth}
\includegraphics[height=4cm]{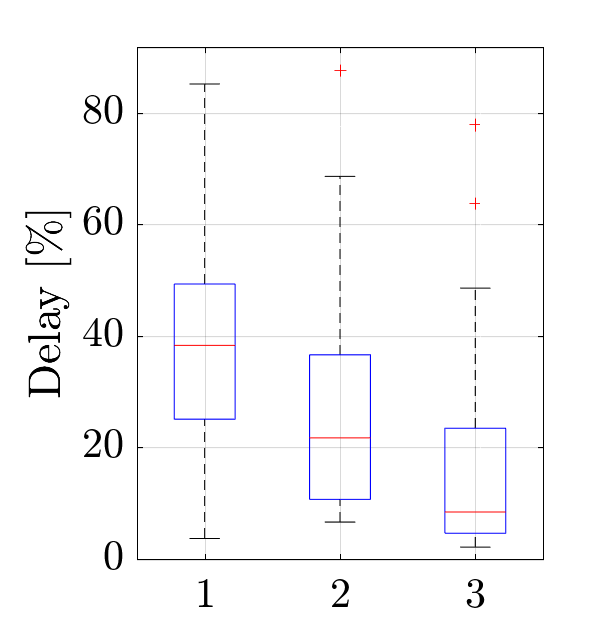}
    \caption{All classes}
    \label{fig:box_ctrlideal}
\end{subfigure}
\caption[]{Box plots showing the increase in TTS compared to the ideal actuation case.}
\label{fig:boxplots}

\begin{center}
\captionof{table}{Average and median TTS for each vehicle class and all vehicles.}
\label{tab:TTS}
\begin{tabular}{|c||c|c||c|c||c|c||c|c|}
\hline
TTS & \multicolumn{ 2}{c||}{Class $a$} & \multicolumn{ 2}{c||}{Class $b$} & \multicolumn{ 2}{c||}{Class $c$} & \multicolumn{ 2}{c|}{Total} \\
\cline{2-9}
[veh $h$] & average & median & average & median & average & median & average & median\\
\hline
Case $(a)$ & $22.62$ & $22.94$ & $369.84$ & $374.13$ & $56.62$ & $56.04$ & $449.08$ & $453.94$ \\    
\hline
Case $(b)$ & $23.25$ & $23.03$ & $329.91$ & $315.75$ &  $60.62$ & $60.37$ & $413.78$ & $398.18$\\  
\hline
Case $(c)$ & $21.77$ & $21.34$ & $304.90$ & $278.63$  & $58.60$ & $58.41$  & $385.27$ & $357.58$ \\ 
\hline
Case $(d)$ & $17.00$ & $16.91$ & $255.00$ & $254.09$ & $55.92$ & $55.93$ &  $327.92$ & $326.42$\\ 
\hline
\end{tabular}
\end{center}


\begin{center}
\captionof{table}{Average and median delay for individual vehicle classes and all vehicles.}
\label{tab:delay}
\begin{tabular}{|c||c|c||c|c||c|c||c|c|}
\hline
Delay & \multicolumn{ 2}{c||}{Class $a$} & \multicolumn{ 2}{c||}{Class $b$} & \multicolumn{ 2}{c||}{Class $c$} & \multicolumn{ 2}{c|}{Total} \\
\cline{2-9}
[$\%$] & average & median & average & median & average & median & average & median\\
\hline
Case $(a)$ & $33.1$ & $35.3$ & $45.0$ & $46.9$ & $1.2 $ & $0.0$ & $36.9$ & $38.3$  \\
\hline
Case $(b)$ & $36.8$ & $33.5$ & $29.4$ & $24.1$ & $8.4$ & $8.0$ & $26.2$ & $21.7$\\ 
\hline
Case $(c)$ & $28.1$ & $20.7$ & $19.6$ & $8.6$  & $4.8$ & $4.4$  & $17.5$ & $8.4$ \\ 
\hline
\end{tabular}
\end{center}

\end{figure*}

Consider the uncontrolled case shown in Figure~\ref{fig:ex_noctrl}.
Around time $t=0.144$ h, the aggregate density of the platooned vehicles and background traffic arriving at the bottleneck is too high, and the aggregate demand exceeds bottleneck capacity.
This causes a traffic breakdown, and after a brief transient, congestion is formed and bottleneck capacity is reduced.
Because of this, even though the incoming traffic density is lower after $t=0.154$ h, and would not exceed the original bottleneck capacity, it is not enough to dissipate the congestion at the bottleneck.
Consequently, the throughput is reduced, the total time spent significantly increased, and the bottleneck will stay congested until the inflow to the highway segment is reduced close to the end of the simulation run.

In contrast to this, in the ideally actuated case shown in Figure~\ref{fig:ex_ctrlideal}, a part of the mainstream-bound background traffic is delayed so that the density of other vehicles reaching the bottleneck while a platoon is traversing it is low enough so as not to cause traffic breakdown and capacity drop.
In this way, free flow is maintained and throughput is close to its theoretical maximum.

As shown in Figure~\ref{fig:ex_ctrlnoramp} and Figure~\ref{fig:ex_ctrlwramp}, the performance of the two proposed control laws is similar.
However, in case the influence of on- and off-ramps is ignored while predicting the evolution of the system, the applied control action is more severe than required, which leads to more congestion upstream of the off-ramp and overall lower efficiency.
The control law that takes the on- and off-ramps into account comes close to emulating the ideal actuation case, but achieves somewhat worse performance because it is unable to selectively affect only one class of background traffic, only has access to the average splitting ratio for the off-ramp, and requires delaying the platoons.


We executed $50$ Monte Carlo simulations, with the same platoon arrival times and background traffic inflow profiles for each control case.
The resulting average and median TTS are shown in Table~\ref{tab:TTS}.
We show the TTS of each vehicle class, and for all vehicles combined.
Apart from comparing the TTS, we also considered the delay, defined as the difference in TTS compared to the ideal actuation case, which is taken as a benchmark for minimum achievable TTS of each simulation run.
The delay is shown as percentage of minimum travel time, and it is shown as a box plot in Figure~\ref{fig:boxplots}, and given in Table~\ref{tab:delay}.
For example, if a vehicle would traverse the road segment in $3$ minutes if it travelled at free flow speed, and actually traverses it in $4.5$ minutes, we say that it had a  $50\%$ delay.

We can see that even by applying control that ignores the existence of on- and off-ramps, as described in Section~\ref{sec:ctrl_noramps}, we reduce the TTS by about $10\%$ of the ideal TTS on average, with the median reduced by about $17\%$.
This corresponds to eliminating $29.1\%$ of the delay on average, or $43.7\%$ by median.
However, only the TTS of class $b$, the mainstream-bound background traffic, is reduced, while the TTS of other vehicles is even somewhat increased.
This can be explained by the fact that the controller assumes that all vehicles are headed for the bottleneck, and will therefore delay the traffic too much, stalling the off-ramp-bound traffic which would otherwise be able to leave the highway unhindered.
In spite of this inefficiency, and owing to the fact that vehicles of class $b$ comprise the majority of the traffic, this control law is still able to preserve free flow and forestall capacity drop at the bottleneck, thus the overall TTS and delays are lower than in the uncontrolled case.

In contrast, when the control from Section~\ref{sec:ctrl_wramps} was used, the TTS of both class $a$ (the platooned vehicles) and class $b$ vehicles, is reduced, with the aggregate TTS lower by almost $20\%$ of ideal TTS on average, or by almost $30\%$ in median.
This corresponds to eliminating $52.7\%$ of the delay on average, or $75.6\%$ by median.
Even though the platoons will be delayed in order to actuate the control action, their TTS will be lower, since they will avoid waiting in congestion upstream of the bottleneck.
This is especially important, since it shows that it is beneficial for the platooned vehicles to employ this control law, even if their goal is not to optimize the overall traffic performance, but to minimize only their own travel time.
The TTS of class $c$ vehicles is still increased compared to the uncontrolled case, but less so than with the previous control law.
Overall, this control law comes very close to the ideal case, with the median delay being only $8.4\%$, and an average delay of $17.5\%$.

It is notable that while the proposed control laws achieve significant reduction of both average and median TTS, but there is a number of outliers corresponding to particularly unfavourable simulation runs.
Since the arrival of platoons is modelled as a Poisson process, we can expect to occasionally have long gaps between two platoons.
If this occurrence coincides with a higher demand of mainstream-bound background traffic, we will not be able to prevent the traffic breakdown, since there would be no platoons available to actuate the control action, resulting in a build-up of congestion and higher a TTS.

\section{Conclusion}
\label{sec:conclusion}

In this work, we used the multi-class CTM framework to study the effects of platoons arriving at a highway bottleneck.
We proposed a simplified queuing model that captures the important aspects of traffic dynamics, and used it to design control laws that use platoons as actuators, and their speed and depth as control inputs, to keep the bottleneck in free flow, maximizing throughput and minimizing total time spent of all vehicles.
The performance of these control laws was tested in multi-class CTM simulations, on a $5$~km long stretch of highway upstream of the lane drop bottleneck, going from three lanes to two.
The considered highway segment also included an on-ramp and an off-ramp.
The achieved TTS using these control laws was compared to the case when no control is used, as well as with the case when we have ideal actuation, and can fully control all individual vehicles.
It has been demonstrated that applying the proposed control laws significantly reduces the TTS compared to the situation with no control, coming close to the performance of the ideal actuation case.
Moreover, even the platooned vehicles, which are delayed in order to affect the rest of traffic, incur lower delays, since they avoid having to traverse the congestion at the bottleneck, making the proposed control beneficial for all traffic participants.

For future work, we are interested in deriving theoretical bounds on effects of the proposed control laws in terms of total time spent and further refining the bound on the achievable throughput.
Additionally, we plan to extend this work to handle longer highway sections, where multiple bottlenecks need to be regulated in cascade, as well as test the approach in microscopic traffic simulations.
In general, the influence of truck platoons on the rest of traffic needs to be further investigated using both simulations and experiments on public roads.


\bibliographystyle{IEEEtran}

{\footnotesize
\bibliography{../traffic_refs, ../disable_url}

\begin{thebibliography}{10}
\providecommand{\url}[1]{#1}
\csname url@samestyle\endcsname
\providecommand{\newblock}{\relax}
\providecommand{\bibinfo}[2]{#2}
\providecommand{\BIBentrySTDinterwordspacing}{\spaceskip=0pt\relax}
\providecommand{\BIBentryALTinterwordstretchfactor}{4}
\providecommand{\BIBentryALTinterwordspacing}{\spaceskip=\fontdimen2\font plus
\BIBentryALTinterwordstretchfactor\fontdimen3\font minus
  \fontdimen4\font\relax}
\providecommand{\BIBforeignlanguage}[2]{{%
\expandafter\ifx\csname l@#1\endcsname\relax
\typeout{** WARNING: IEEEtran.bst: No hyphenation pattern has been}%
\typeout{** loaded for the language `#1'. Using the pattern for}%
\typeout{** the default language instead.}%
\else
\language=\csname l@#1\endcsname
\fi
#2}}
\providecommand{\BIBdecl}{\relax}
\BIBdecl

\bibitem{bergenhem2012overview}
C.~Bergenhem, S.~Shladover, E.~Coelingh, C.~Englund, and S.~Tsugawa, ``Overview
  of platooning systems,'' in \emph{Proceedings of the 19th ITS World
  Congress}, Vienna, Austria, 2012.

\bibitem{bonnet2000fuel}
C.~Bonnet and H.~Fritz, ``Fuel consumption reduction in a platoon: Experimental
  results with two electronically coupled trucks at close spacing,'' SAE
  Technical Paper, Tech. Rep., 2000.

\bibitem{mckinsey2018platooning}
\BIBentryALTinterwordspacing
A.~Chottani, G.~Hastings, J.~Murnane, and F.~Neuhaus, ``Distraction or
  disruption? autonomous trucks gain ground in us logistics,'' \emph{{\relax
  McKinsey Smith Co.}}, 2018. [Online]. Available:
  \url{https://www.mckinsey.com/industries/travel-transport-and-logistics/our-insights/distraction-or-disruption-autonomous-trucks-gainground-in-us-logistics}
\BIBentrySTDinterwordspacing

\bibitem{shladover2012impacts}
S.~Shladover, D.~Su, and X.-Y. Lu, ``Impacts of cooperative adaptive cruise
  control on freeway traffic flow,'' \emph{Transportation Research Record:
  Journal of the Transportation Research Board}, no. 2324, pp. 63--70, 2012.

\bibitem{lioris2017platoons}
J.~Lioris, R.~Pedarsani, F.~Y. Tascikaraoglu, and P.~Varaiya, ``Platoons of
  connected vehicles can double throughput in urban roads,''
  \emph{Transportation Research Part C: Emerging Technologies}, vol.~77, pp.
  292--305, 2017.

\bibitem{moridpour2015impact}
S.~Moridpour, E.~Mazloumi, and M.~Mesbah, ``Impact of heavy vehicles on
  surrounding traffic characteristics,'' \emph{Journal of advanced
  transportation}, vol.~49, no.~4, pp. 535--552, 2015.

\bibitem{alam2010experimental}
A.~Alam, A.~Gattami, and K.~H. Johansson, ``An experimental study on the fuel
  reduction potential of heavy duty vehicle platooning,'' in \emph{13th
  International IEEE Conference on Intelligent Transportation Systems}.\hskip
  1em plus 0.5em minus 0.4em\relax IEEE, 2010, pp. 306--311.

\bibitem{aarts2016european}
L.~Aarts and G.~Feddes, ``European truck platooning challenge,'' in
  \emph{HVTT14: International Symposium on Heavy Vehicle Transport Technology,
  Rotorua, New Zealand}, 2016.

\bibitem{bhoopalam2018planning}
A.~K. Bhoopalam, N.~Agatz, and R.~Zuidwijk, ``Planning of truck platoons: A
  literature review and directions for future research,'' \emph{Transportation
  Research Part B: Methodological}, vol. 107, pp. 212--228, 2018.

\bibitem{rijkswaterstaat2016platooning}
{\relax Rijkswaterstaat}, ``European truck platooning challenge 2016: Lessons
  learnt,'' 2016.

\bibitem{duret2019hierarchical}
A.~Duret, M.~Wang, and A.~Ladino, ``A hierarchical approach for splitting truck
  platoons near network discontinuities,'' \emph{Transportation Research Part
  B: Methodological}, 2019.

\bibitem{jin2018modeling}
L.~Jin, M.~{\v{C}}i{\v{c}}i{\v{c}}, S.~Amin, and K.~H. Johansson, ``Modeling
  the impact of vehicle platooning on highway congestion: A fluid queuing
  approach,'' in \emph{Proceedings of the 21st International Conference on
  Hybrid Systems: Computation and Control (part of CPS Week)}.\hskip 1em plus
  0.5em minus 0.4em\relax ACM, 2018, pp. 237--246.

\bibitem{wang2014local}
Y.~Wang, E.~B. Kosmatopoulos, M.~Papageorgiou, and I.~Papamichail, ``Local ramp
  metering in the presence of a distant downstream bottleneck: Theoretical
  analysis and simulation study,'' \emph{IEEE Transactions on Intelligent
  Transportation Systems}, vol.~15, no.~5, pp. 2024--2039, 2014.

\bibitem{hadiuzzaman2012variable}
M.~Hadiuzzaman, T.~Z. Qiu, and X.-Y. Lu, ``Variable speed limit control design
  for relieving congestion caused by active bottlenecks,'' \emph{Journal of
  Transportation Engineering}, vol. 139, no.~4, pp. 358--370, 2012.

\bibitem{herrera2010evaluation}
J.~C. Herrera, D.~B. Work, R.~Herring, X.~J. Ban, Q.~Jacobson, and A.~M. Bayen,
  ``Evaluation of traffic data obtained via {GPS}-enabled mobile phones: The
  mobile century field experiment,'' \emph{Transportation Research Part C:
  Emerging Technologies}, vol.~18, no.~4, pp. 568--583, 2010.

\bibitem{cui2017stabilizing}
S.~Cui, B.~Seibold, R.~Stern, and D.~B. Work, ``Stabilizing traffic flow via a
  single autonomous vehicle: Possibilities and limitations,'' in
  \emph{Intelligent Vehicles Symposium (IV)}.\hskip 1em plus 0.5em minus
  0.4em\relax IEEE, 2017, pp. 1336--1341.

\bibitem{wu2018stabilizing}
C.~Wu, A.~M. Bayen, and A.~Mehta, ``Stabilizing traffic with autonomous
  vehicles,'' in \emph{2018 IEEE International Conference on Robotics and
  Automation (ICRA)}.\hskip 1em plus 0.5em minus 0.4em\relax IEEE, 2018, pp.
  1--7.

\bibitem{piacentini2018traffic}
G.~Piacentini, P.~Goatin, and A.~Ferrara, ``Traffic control via moving
  bottleneck of coordinated vehicles,'' in \emph{15th IFAC symposium on control
  in transportation systems}, 2018.

\bibitem{vinitsky2018lagrangian}
E.~Vinitsky, K.~Parvatey, A.~Kreidiehz, C.~Wu, and A.~Bayen, ``Lagrangian
  control through deep-rl: Applications to bottleneck decongestion,'' in
  \emph{21st IEEE International Conference on Intelligent Transportation
  Systems (ITSC)}, Maui, US, 2018.

\bibitem{cicic2018traffic}
M.~{{\v{C}}i{\v{c}}i{\'c}} and K.~H. Johansson, ``Traffic regulation via
  individually controlled automated vehicles: a cell transmission model
  approach,'' in \emph{21st IEEE International Conference on Intelligent
  Transportation Systems}, Maui, US, 2018.

\bibitem{alam2015heavy}
A.~Alam, B.~Besselink, V.~Turri, J.~Martensson, and K.~H. Johansson,
  ``Heavy-duty vehicle platooning for sustainable freight transportation: A
  cooperative method to enhance safety and efficiency,'' \emph{IEEE Control
  Systems}, vol.~35, no.~6, pp. 34--56, 2015.

\bibitem{liu2016modeling}
L.~Liu, L.~Zhu, and D.~Yang, ``Modeling and simulation of the car-truck
  heterogeneous traffic flow based on a nonlinear car-following model,''
  \emph{Applied Mathematics and Computation}, vol. 273, pp. 706--717, 2016.

\bibitem{delle2014scalar}
M.~L. Delle~Monache and P.~Goatin, ``Scalar conservation laws with moving
  constraints arising in traffic flow modeling: an existence result,''
  \emph{Journal of Differential equations}, vol. 257, no.~11, pp. 4015--4029,
  2014.

\bibitem{cicic2019multiclass}
M.~{\v{C}}i{\v{c}}i{\'c} and K.~H. Johansson, ``Stop-and-go wave dissipation
  using accumulated controlled moving bottlenecks in multi-class ctm
  framework,'' in \emph{Control and Decision Conference (ECC)}, Nice, France,
  2019.

\bibitem{daganzo1994cell}
C.~F. Daganzo, ``The cell transmission model: A dynamic representation of
  highway traffic consistent with the hydrodynamic theory,''
  \emph{Transportation Research Part B: Methodological}, vol.~28, no.~4, pp.
  269--287, 1994.

\bibitem{piacentini2019mcctm}
G.~Piacentini, M.~\v{C}i\v{c}i\'{c}, A.~Ferrara, and K.~H. Johansson, ``{VACS}
  equipped vehicles for congestion dissipation in multi-class {CTM}
  framework,'' in \emph{European Control Conference}, 2019.

\bibitem{kontorinaki2017first}
M.~Kontorinaki, A.~Spiliopoulou, C.~Roncoli, and M.~Papageorgiou, ``First-order
  traffic flow models incorporating capacity drop: Overview and real-data
  validation,'' \emph{Transportation Research Part B: Methodological}, vol.
  106, pp. 52--75, 2017.

\bibitem{han2016new}
Y.~Han, Y.~Yuan, A.~Hegyi, and S.~P. Hoogendoorn, ``New extended discrete
  first-order model to reproduce propagation of jam waves,''
  \emph{Transportation Research Record: Journal of the Transportation Research
  Board}, no. 2560, pp. 108--118, 2016.

\end{thebibliography}
}

\begin{IEEEbiography}[{\includegraphics[width=1in,height=1.25in,clip,keepaspectratio]{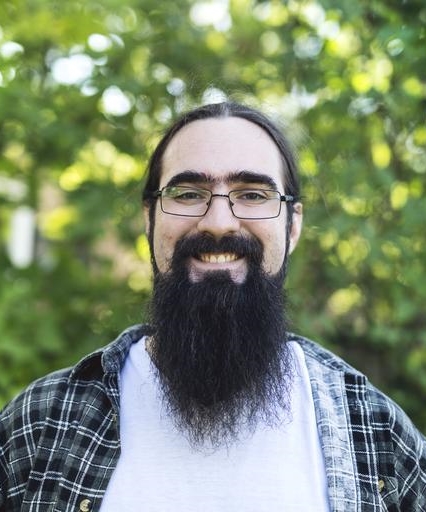}}]{Mladen \v Ci\v ci\'c} received his B.Sc. and M.Sc. degrees in electrical engineering and computer science, automatic control, from the School of Electrical Engineering, University of Belgrade. He is pursuing the Ph.D degree at the Division of Decision and Control Systems, KTH Royal Institute of Technology. He was a member of Marie Sk{\l}odowska Curie oCPS ITN and an affiliated Wallenberg AI, Autonomous Systems and Software Program (WASP) student. He was a visiting scholar at the C2SMART University Transportation Center in the New York University Tandon School of Engineering. His research interests include traffic control using mixed traffic models. 
\end{IEEEbiography}

\begin{IEEEbiography}[{\includegraphics[width=1in,height=1.25in,clip,keepaspectratio]{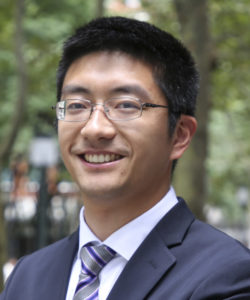}}]{Li Jin}
is an Assistant Professor at the Department of Civil and Urban Engineering and the C2SMART University Transportation Center in the New York University Tandon School of Engineering.
He received B.Eng. degree in Mechanical Engineering from Shanghai Jiao Tong University in 2011, M.S. degree in Mechanical Engineering from Purdue University in 2012, and Ph.D. degree in Transportation from the Massachusetts Institute of Technology in 2018.
His research focuses on developing resilient control algorithms for cyber-physical systems with guarantees of efficiency in nominal settings, robustness against random perturbations, and survivability under strategic disruptions.
He is a recipient of the Ho-Ching and Hang-Ching Fund Award and Schoettler Scholarship Fund.
\end{IEEEbiography}

\begin{IEEEbiography}[{\includegraphics[width=1in,height=1.25in,clip,keepaspectratio]{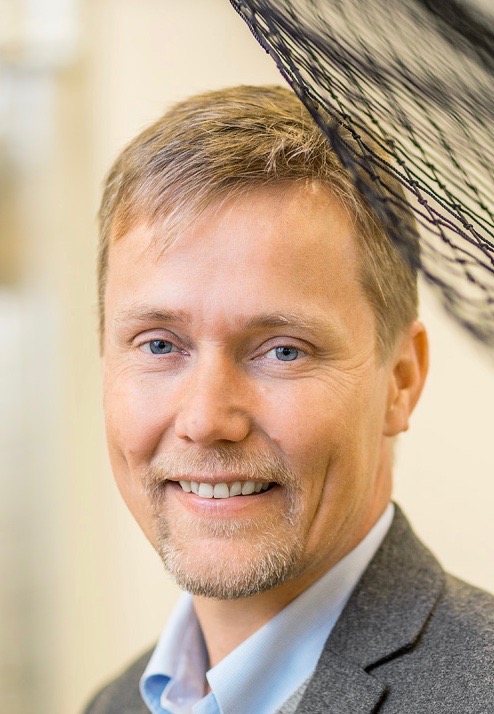}}]{Karl H. Johansson} is Director of the
Stockholm Strategic Research Area ICT -- the Next
Generation and Professor at the School of Electrical
Engineering, KTH Royal Institute of Technology.
He received MSc and PhD degrees in Electrical
Engineering from Lund University. He has held
visiting positions at UC Berkeley, Caltech, NTU,
HKUST Institute of Advanced Studies, and NTNU.
His research interests are in networked control
systems, cyber-physical systems, and applications
in transportation, energy, and automation. He is a
member of the IEEE Control Systems Society Board of Governors and the
European Control Association Council. He has received several best paper
awards and other distinctions, including a ten-year Wallenberg Scholar Grant,
a Senior Researcher Position with the Swedish Research Council, the Future
Research Leader Award from the Swedish Foundation for Strategic Research,
and the triennial Young Author Prize from IFAC. He is member of the Royal
Swedish Academy of Engineering Sciences, Fellow of the IEEE, and IEEE
Distinguished Lecturer.
\end{IEEEbiography}

\end{document}